\documentclass[leqno,11pt]{article}
\usepackage{a4wide,amsmath,amssymb,amsfonts,mathrsfs,exscale,graphics,amsthm,colonequals,comment}

\usepackage[curve,matrix,arrow,cmtip]{xy}
\newdir^{ (}{{}*!/-3pt/\dir^{(}}   
\newdir_{ (}{{}*!/-3pt/\dir_{(}}   
\NoComputerModernTips

\newcommand{\Ghat}{\hat G}

\newcommand{\Gbar}{\bar G}

\numberwithin{equation}{section}
\def\UseTheoremCounterForNextEquation{\setcounter{equation}{\value{theorem}}\addtocounter{theorem}{1}}

\theoremstyle{plain}
\newtheorem{theorem}{Theorem}[section]
\newtheorem{lemma}[theorem]{Lemma}
\newtheorem{corollary}[theorem]{Corollary}
\newtheorem{proposition}[theorem]{Proposition}

\theoremstyle{definition}
\newtheorem{definition}[theorem]{Definition}

\newtheorem{construction}[theorem]{Construction}
\newtheorem{remark}[theorem]{Remark}
\newtheorem{example}[theorem]{Example}

\newcommand\lto{\longrightarrow}
\newcommand\ltoover[1]{\mathrel{\smash{\overset{#1}{\lto}}}}
\newcommand\varto[1]{\mathrel{\hbox to #1pt{\rightarrowfill}}}

\newcommand{\lleftto}{\longleftarrow}
\newcommand\llefttoover[1]{\mathrel{\smash{\overset{#1}{\lleftto}}}}

\newcommand{\sends}{\mapsto}

\newcommand{\iso}{\overset{\sim}{\to}}
\newcommand{\liso}{\overset{\sim}{\lto}}

\renewcommand{\implies}{\Rightarrow}

\let\longto\longrightarrow
\let\onto\twoheadrightarrow
\def\isoto{\stackrel{\sim}{\longto}}

\let\phi\varphi
\let\epsilon\varepsilon
\let\setminus\smallsetminus
\let\emptyset\varnothing

\let\ge\geqslant

\def\Norm{\mathop{\rm Norm}\nolimits}

\def\Stab{\mathop{\rm Stab}\nolimits}

\def\Lie{\mathop{\rm Lie}\nolimits}
\def\Ad{\mathop{\rm Ad}\nolimits}

\def\GL{\mathop{\rm GL}\nolimits}

\def\codim{\mathop{\rm codim}\nolimits}

\def\id{{\rm id}}

\newcommand{\BF}{{\mathbb{F}}}

\newcommand{\CA}{{\mathcal A}}

\newcommand{\CR}{{\mathcal R}}

\newcommand{\CZ}{{\mathcal Z}}

\newcommand{\ru}{\mathcal{R}_u}
\newcommand{\defeq}{\colonequals}

\newcommand{\leftexp}[2]{{\vphantom{#2}}^{#1}{#2}}
\newcommand{\doubleexp}[3]{\leftexp{#1}{#2}^{#3}}
\newcommand{\doubleexpIJ}[1]{\,{\vphantom{#1}}^{I}\!{#1}^{J}}

\theoremstyle{definition}

\DeclareMathOperator{\inn}{int}
\DeclareMathOperator{\orb}{o}
\DeclareMathOperator{\relpos}{relpos}
\DeclareMathOperator{\Par}{Par}

\newcommand{\BIGOP}[1]{\mathop{\mathchoice%
{\raise-0.22em\hbox{\huge $#1$}}%
{\raise-0.05em\hbox{\Large $#1$}}{\hbox{\large $#1$}}{#1}}}

\newcommand{\BIGboxplus}{\mathop{\mathchoice%
{\raise-0.35em\hbox{\huge $\boxplus$}}%
{\raise-0.15em\hbox{\Large $\boxplus$}}{\hbox{\large $\boxplus$}}{\boxplus}}}

\newcounter{listcounter}
\newcounter{deflistcounter}
\newcounter{equivcounter}

\newskip{\itemsepamount}
\newskip{\topsepamount}
% ...

\begin{document}
\newcounter{ctr} 

\title{Algebraic zip data}

\author{%
\renewcommand{\thefootnote}{\arabic{footnote}}%
Richard Pink\footnote{Dept. of Mathematics,
 ETH Z\"urich, 
 CH-8092 Z\"urich,
 Switzerland,
 \tt pink@math.ethz.ch} \hspace{1cm}
 Torsten Wedhorn\footnote{Dept. of Mathematics,
 University of Paderborn, 
 D-33098 Paderborn,
 Germany,
 \tt wedhorn@math.uni-paderborn.de} \hspace{1cm}
 Paul Ziegler\footnote{Dept. of Mathematics,
 ETH Z\"urich
 CH-8092 Z\"urich,
 Switzerland,
 \tt paul.ziegler@math.ethz.ch}
}

\date{\today}

\maketitle


\begin{abstract}
An algebraic zip datum is a tuple $\CZ = (G,P,Q,\varphi)$ consisting of a reductive group $G$ together with parabolic subgroups $P$ and $Q$ and an isogeny $\varphi\colon P/R_uP\to Q/R_uQ$. We study the action of the group $E_\CZ := \bigl\{ (p,q)\in P{\times}Q \bigm| \varphi(\pi_{P}(p)) =\pi_Q(q)\bigr\}$ on $G$ given by $((p,q),g)\mapsto pgq^{-1}$. 
We define certain smooth $E_\CZ$-invariant subvarieties of $G$, show that they define a stratification of $G$. We determine their dimensions and their closures and give a description of the stabilizers of the $E_\CZ$-action on $G$. We also generalize all results to non-connected groups.

We show that for special choices of $\CZ$ the algebraic quotient stack $[E_\CZ \backslash G]$ is isomorphic to $[G \backslash Z]$ or to $[G \backslash Z']$, where $Z$ is a $G$-variety studied by Lusztig and He in the theory of character sheaves on spherical compactifications of $G$ and where $Z'$ has been defined by Moonen and the second author in their classification of $F$-zips. In these cases the $E_\CZ$-invariant subvarieties correspond to the so-called ``$G$-stable pieces'' of $Z$ defined by Lusztig (resp.~the $G$-orbits of $Z'$). 
\end{abstract}


\section{Introduction}\label{Intro}


\subsection{Background}
\label{background}
Let $G$ be a connected reductive linear algebraic group over an algebraically closed field~$k$. Then $G \times G$ acts on $G$ via simultaneous left and right translation $((g_1,g_2),g)\mapsto g_1gg_2^{-1}$. In a series of papers, Lusztig (\cite{Lusztig:parsheaves1}, \cite{Lusztig:parsheaves2}), He (\cite{He:GStablePieces}, \cite{He:CharSheavesGroupComp}, \cite{He:MinimalLengthCoxeter}), and Springer (\cite{Springer:BruhatLemma}) studied a certain spherical $G \times G$-equivariant smooth compactification $\Gbar$ of~$G$. For $G$ semi-simple adjoint this is the so-called wonderful compactification from~\cite{conciniprocesi}. In general the $G \times G$-orbits $Z_I\subset\Gbar$ are in natural bijection to the subsets $I$ of the set of simple reflections in the Weyl group of~$G$. Lusztig and He defined and studied so-called $G$-stable pieces in~$Z_I$, which are certain subvarieties that are invariant under the diagonally embedded subgroup $G\subset G \times G$. These $G$-stable pieces play an important role in their study of character sheaves on~$\Gbar$. Lusztig and He also consider non-connected groups, corresponding to twisted group actions. Other generalizations of these varieties have been considered by Lu and Yakimow (\cite{LuYakimov}). A further motivation to study $G$-stable pieces comes from Poisson geometry: It was proved by Evens and Lu (\cite{EvensLu:LagrangianII}), that for certain Poisson structure, each $G$-orbit on $Z_I$ is a Poisson submanifold.

In~\cite{moonwed} Moonen and the second author studied the De Rham cohomology $H^{\bullet}_{\rm DR}(X/k)$ of a smooth proper scheme $X$ with degenerating Hodge spectral sequence over an algebraically closed field $k$ of positive characteristic. They showed that $H^{\bullet}_{\rm DR}(X/k)$ carries the structure of a so-called \emph{$F$-zip}, namely: it is a finite-dimensional $k$-vector space together with two filtrations (the ``Hodge'' and the ``conjugate'' filtration) and a Frobenius linear isomorphism between the associated graded vector spaces (the ``Cartier isomorphism''). They showed that the isomorphism classes of $F$-zips of fixed dimension $n$ and with fixed type of Hodge filtration are in natural bijection with the orbits under $G\defeq\GL_{n,k}$ of a variant $Z'_I$ of the $G \times G$-orbit $Z_I$ studied by Lusztig. They studied the varieties $Z_I'$ for arbitrary reductive groups $G$ and determined the $G$-orbits in them as analogues of the $G$-stable pieces in~$Z_I$. By specializing $G$ to classical groups they deduce from this a classification of $F$-zips with additional structure, e.g., with a non-degenerate symmetric or alternating form. They also consider non-connected groups. 
Moreover, the automorphism group of an $F$-zip is isomorphic to the stabilizer in $G$ of any corresponding point in~$Z'_I$. 

When $X$ varies in a smooth family over a base~$S$, its relative De Rham cohomology forms a family of $F$-zips over~$S$. The set of points of $S$ where the $F$-zip lies in a given isomorphism class is a natural generalization of an Ekedahl-Oort stratum in the Siegel moduli space. Information about the closure of such a stratum corresponds to information about how the isomorphism class of an $F$-zip can vary in a family, and that in turn is equivalent to determining which $G$-orbits in $Z_I'$ are contained in the closure of a given $G$-orbit. 

In each of these cases one is interested in the classification of the $G$-stable pieces, the description of their closures, and the stabilizers of points in~$G$. In this article we give a uniform approach to these questions that generalizes all the above situations. 


\subsection{Main results}
\label{MainResults}

The central definition in this article is the following:

\begin{definition}\label{1ZipDatumDef}
A \emph{connected algebraic zip datum} is a tuple $\CZ = (G,P,Q,\varphi)$ consisting of a connected reductive linear algebraic group $G$ over~$k$ together with parabolic subgroups $P$ and $Q$ and an isogeny $\varphi\colon P/R_uP\to Q/R_uQ$. The group
$$E_\CZ := \bigl\{ (p,q)\in P{\times}Q \bigm| \phi(\pi_{P}(p)) =\pi_Q(q)\bigr\}$$
is called the \emph{zip group associated to~$\CZ$}. It acts on $G$ through the map $((p,q),g)\mapsto pgq^{-1}$.
The union of the $E_\CZ$-orbits of all elements of a subset $X \subset G$ is denoted by $\orb_\CZ(X)$.
\end{definition}

Fix such data $\CZ = (G,P,Q,\varphi)$. To apply the machinery of Weyl groups to $\CZ$ we choose a Borel subgroup $B$ of~$G$, a maximal torus $T$ of~$B$, and an element $g$ of $G$ such that $B\subset Q$, $\leftexp{g}{B} \subset P$, $\phi\bigl(\pi_P(\leftexp{g}{B})\bigr) = \pi_Q(B)$, and $\phi\bigl(\pi_P(\leftexp{g}{T})\bigr) = \pi_Q(T)$. Let $W$ denote the Weyl group of $G$ with respect to~$T$, and $S\subset W$ the set of simple reflections corresponding to~$B$. Let $I\subset S$ be the type of the parabolic~$P$ and $W_I\subset W$ its Weyl group. Let $\leftexp{I}{W}$ be the set of all $w\in W$ that have minimal length in their coset~$W_Iw$. To each $w\in\leftexp{I}{W}$ we associate the $E_{\CZ}$-invariant subset 
\UseTheoremCounterForNextEquation
\begin{equation}\label{1GwDef}
G^{w} = \orb_\CZ(gB\dot wB)
\end{equation}
and prove (Theorems \ref{thm:zstablepieces1}, \ref{GwVarDim} and \ref{Gw4}):

\begin{theorem}\label{1Decomposition}
The $E_\CZ$-invariant subsets $G^w$ form a pairwise disjoint decomposition of $G$ into locally closed smooth subvarieties. The dimension of $G^w$ is $\dim P + \ell(w)$.
\end{theorem}

Next the isogeny $\varphi$ induces an isomorphism of Coxeter system $\psi\colon (W_I,I) \iso (W_J,J)$ (see~\eqref{DefinePsi} for its precise definition), where $W_J \subset W$ and $J \subset S$ are the Weyl group and the type of the parablic subgroup~$Q$. Let $\leq$ denote the Bruhat order on~$W$. We prove (Theorem~\ref{thm:closure}):

\begin{theorem}\label{1Closure}
The closure of $G^w$ is the union of $G^{w'}$ for all $w'\in \leftexp{I}{W}$ such that there exists $y \in W_I$ with $yw'\psi(y)^{-1} \leq w$.
\end{theorem}

We call $\CZ$ \emph{orbitally finite} if the number of $E_\CZ$-orbits in $G$ is finite. We give a necessary and sufficient criterion for this to happen (Proposition~\ref{OrbFinEqu}). In particular it happens when the differential of $\varphi$ at $1$ vanishes, for instance if $\varphi$ is a Frobenius isogeny (Proposition~\ref{dphi=0}). We prove (Theorem~\ref{thm:representatives}):

\begin{theorem}\label{1representatives}
If $\CZ$ is orbitally finite, then each $G^w$ is a single $E_\CZ$-orbit, and so the set $\{g\dot w\mid w\in \leftexp{I}{W}\}$ is a set of representatives for the $E_\CZ$-orbits in $G$.
\end{theorem}

One can also consider the $E_\CZ$-orbit of $g\dot w$ for any element $w\in W$ instead of just those in~$\leftexp{I}{W}$. It is then natural to ask when two such orbits lie in the same $E_\CZ$-invariant piece~$G^w$. (For orbitally finite $\CZ$ this is equivalent to asking when the orbits are equal.) We give an explicit description of this equivalence relation on $W$ that depends only on the subgroup $W_I$ and the homomorphism~$\psi$ (Theorem~\ref{gut}). We prove that all equivalence classes have the same cardinality~$\#W_I$, although they are in general no cosets of~$W_I$ and we do not know a simple combinatorial description for them. It is intriguing that we obtain analogous results for an \emph{abstract zip datum} based on an arbitrary finitely generated Coxeter group (Theorem~\ref{Abszstablepieces1}) or even an arbitrary abstract group (Theorem~\ref{eqclasscard}) in place of~$W$.

Other results include information on point stabilizers and infinitesimal stabilizers (Section~\ref{Stabilizer}), the generalization of the main results to non-connected groups (Section~\ref{nonconnected}), a dual parametrization by a set $W^J$ in place of~$\leftexp{I}{W}$ (Section~\ref{Diff}) and the relation with the varieties $Z_I$ studied by Lusztig and He and their generalizations~$Z_I'$ (Section~\ref{CompareLusztig}).


\subsection{Applications}
\label{Applications}

Let us explain why this theory of algebraic zip data is a generalization of the situations described in Subsection~\ref{background}. In Section \ref{CompareLusztig} we consider a connected reductive algebraic group $G$ over~$k$, an isogeny $\phi\colon G \to G$, a subset $I$ of the set of simple reflections associated to~$G$, and an element $x$ in the Weyl group of $G$ satisfying certain technical conditions. To such data we associate a certain algebraic variety $X_{I,\phi,x}$ with an action of~$G$, a certain connected algebraic zip datum $\CZ$ with underlying group~$G$, and morphisms
\UseTheoremCounterForNextEquation
\begin{equation}\label{EqPsi}
G \llefttoover{\lambda} G\times G \ltoover{\rho} X_{I,\phi,x}
\end{equation}

In Theorem~\ref{orbitbijection} we show that there is a closure preserving bijection between the $E_\CZ$-invariant subsets of $A\subset G$ and the $G$-invariant subsets of $B\subset X_{I,\phi,x}$ given by $\lambda^{-1}(A)=\rho^{-1}(B)$. We also prove that the stabilizer in $E_\CZ$ of $g\in G$ is isomorphic to the stabilizer in $G$ of any point of the $G$-orbits in $X_{I,\phi,x}$ corresponding to the orbit of $g$. These results can also be phrased in the language of algebraic stacks, see Theorem \ref{stacks}.

\medskip
In the special case $\phi = \id_G$ the above $X_{I,\phi,x}$ is the variety $Z_I$ defined by Lusztig. In Theorem~\ref{Lusztigstablepieces} we verify that the subsets $G^w\subset G$ correspond to the $G$-stable pieces defined by him. Thus Theorem \ref{1Closure} translates to a description of the closure relation between these $G$-stable pieces, which had been proved before by He~\cite{He:GStablePieces}.

\medskip
If ${\rm char}(k)$ is positive and $\phi$ is the Frobenius isogeny associated to a model of $G$ over a finite field, the above $X_{I,\phi,x}$ is the variety $Z'_I$ defined by Moonen and the second author. In this case the zip datum $\CZ$ is orbitally finite, and so we obtain the main classification result for the $G$-orbits in~$Z'_I$ from~\cite{moonwed}, the closure relation between these $G$-orbits, and the description of the stabilizers in $G$ of points in~$Z'_I$. In this case the closure relation had been determined in the unpublished note~\cite{specfzips}, the ideas of which are reused in the present article. Meanwhile Viehmann~\cite{Viehmann:Trunc} has given a different proof of the closure relation in this case using the theory of loop groups. For those cases which pertain to the study of modulo $p$ reductions of $F$-crystals with additional structure that show up in the study of special fibers of good integral models of Shimura varieties of Hodge type Moonen (\cite{Moonen:GroupAddStructure}) and, more generally, Vasiu (\cite{Vasiu:ModPClassification}) have obtained similar classification results. In these cases Vasiu (loc.~cit.) has also shown that the connected component of the stabilizers are unipotent.

\medskip
For $G = \GL_n$ (resp.\ a classical group) we therefore obtain a new proof of the classification of $F$-zips (resp.~of $F$-zips with additional structure) from~\cite{moonwed}. We can also deduce how $F$-zips (possibly with additional structure) behave in families, and can describe their automorphism groups as the stabilizers in $E_\CZ$ of the corresponding points of~$G$. This is applied in~\cite{ViWd} to the study of Ekedahl-Oort strata for Shimura varieties of PEL type.


\subsection{Contents of the paper}

In Section~\ref{prelim} we collect some results on algebraic groups and Coxeter systems that are used in the sequel. 
Algebraic zip data $\CZ$ are defined in Section~\ref{AlgZipData}, where we also establish basic properties of the triple $(B,T,g)$, called a \emph{frame of~$\CZ$}. 

Section~\ref{Bruhat} is based on the observation that every $E_\CZ$-orbit is contained in the double coset $Pg\dot xQ$ for some $x\in W$ and meets the subset $g\dot xM$, where $M$ is a Levi subgroup of~$Q$. In it we define another zip datum $\CZ_{\dot x}$ with underlying reductive group~$M$ and establish a number of results relating the $E_\CZ$-orbits in $Pg\dot xQ$ to the $E_{\CZ_{\dot x}}$-orbits in~$M$. This is the main induction step used in most of our results.

In Section~\ref{Gw} we give different descriptions of the $E_\CZ$-invariant subsets $G^w$ for $w\in\leftexp{I}{W}$ and prove Theorem~\ref{1Decomposition}. 
In Section~\ref{Closure} we determine the closure of $G^w$ and prove Theorem~\ref{1Closure}.
Orbitally finite zip data are studied in Section~\ref{Frobenius}, proving Theorem~\ref{1representatives}. 
Section~\ref{Stabilizer} contains some results on point stabilizers and infinitesmial stabilizers. 
Abstract zip data are defined and studied in Section~\ref{Abstract}. In Section~\ref{nonconnected} our main results are generalized to algebraic zip data based on non-connected groups.

In Section \ref{Diff} we discuss a dual parametrization of the subsets $G^w$ by a subset $W^J$ of $W$ in place of~$\leftexp{I}{W}$.
Finally, in Section~\ref{CompareLusztig} we prove the results described in Subsection \ref{Applications} above.

\medskip

The paper is based on parts of the unpublished note \cite{specfzips} by the second author and the master thesis \cite{mt} by the third author, but goes beyond both.
\medskip

After the referee pointed out to us the references \cite{LuYakimov:GroupCompactification} and \cite{He:MinimalLengthCoxeter}, we became aware that there Lu, Yakimov and He study a class of group actions which contains ours when $\phi$ is an isomorphism. In this case, Theorems \ref{1Decomposition} and \ref{1Closure} were already proven in [loc. cit]. Also, many of the ideas we have used to study the decomposition of $G$ into $E_\CZ$-stable pieces are already present there.
\bigskip

\noindent{\scshape Acknowledgements.\ } We thank the referee for pointing out some references. The second author was partially supported by the SPP 1388 ``Representation theory'' of the DFG.


\section{Preliminaries on algebraic groups and Coxeter groups}
\label{prelim}

Throughout, the inner automorphism associated to an element $h$ of a group $G$ will be denoted $\inn(h)\colon {G\to G}$, $g\mapsto \leftexp{h}{g} := h g h^{-1}$. Similarly, for any subset $X\subset G$ we set $\leftexp{h}{X} := hXh^{-1}$. 

\subsection{General facts about linear algebraic groups}\label{LinAlg}

Throughout, we use the language of algebraic varieties over a fixed algebraically closed field~$k$. By an algebraic group $G$ we always mean a linear algebraic group over~$k$. We let $\ru G$ denote the unipotent radical of the identity component of $G$ and $\pi_G\colon G\onto G/\ru G$ the canonical projection. An isogeny between two connected algebraic groups is a surjective homomorphism with finite kernel. 

Consider an algebraic group $G$, an algebraic subgroup $H$ of $G$, and a quasi-projective variety $X$ with a left action of~$H$. Then we denote by $G\times^{H} X$ the quotient of $G \times X$ under the left action of $H$ defined by $h\cdot (g,x)=(gh^{-1},h\cdot x)$, which exists by \cite{efa}, Section 3.2. The action of $G$ on $G\times X$ by left multplication on the first factor induces a left action of $G$ on $G\times^{H} X$. This is the pushout of $X$ with respect to the inclusion $H\hookrightarrow G$.

\begin{lemma}
\label{lem:pushoutorbits}
For $G$, $H$, and $X$ as above, the morphism $X\to G\times^H X$ which sends $x\in X$ to the class of $(1,x)$ induces a closure-preserving bijection between the $H$-invariant subsets of $X$ and the $G$-invariant subsets of $G\times^{H} X$. If $Y\subset X$ is an $H$-invariant subvariety of $X$, then the corresponding $G$-invariant subset of $G\times^H X$ is the subvariety $G\times^H Y$ of $G\times^H X$.
\end{lemma}

\begin{proof}
The morphism in question is the composite of the inclusion $i\colon X\to G\times X$, $x\mapsto (1,x)$ and the projection $\operatorname{pr}\colon G\times X\to G\times^{H} X$. Let $(g,h) \in G\times H$ act on $G\times X$ from the left by $(g^{\prime},x)\mapsto (g g^{\prime}h^{-1},h\cdot x)$. 
Then the $G\times H$-invariant subsets of $G\times X$ are the sets of the form $G\times A$ for $H$-invariant subsets $A\subset X$. Therefore $i$ induces a closure-preserving bijection between the $H$-invariant subsets of $X$ and the $G\times H$-invariant subsets of $G\times X$. Furthermore, since $G\times^{H} X$ carries the quotient topology with respect to $\operatorname{pr}$, the morphism $\operatorname{pr}$ induces a closure-preserving bijection between the $G\times H$-invariant subsets of $G\times X$ and the $G$-invariant of $G\times^{H}X$. Altogether this proves the claim.
\end{proof}

\begin{lemma}[see \cite{slodowy}, Lemma 3.7.4]\label{lem:pushout}
Let $G$ be an algebraic group with an algebraic subgroup~$H$. Let $X$ be a variety with a left action of~$G$. Let $f\colon X\to G/H$ be a $G$-equivariant morphism, and let $Y \subset X$ be the fiber $f^{-1}(H)$. Then $Y$ is stabilized by $H$, and the map $G\times^{H} Y\to X$ sending the equivalence class of $(g,y)$ to $g\cdot y$ defines an isomorphism of $G$-varieties.
\end{lemma}

\begin{lemma}\label{lem:properclosure}
Let $G$ be an algebraic group acting on an algebraic variety $Z$ and let $P\subset G$ be an algebraic subgroup such that $G/P$ is proper. Then for any $P$-invariant subvariety $Y\subset Z$ one has
$$G\cdot \overline Y=\overline{G\cdot Y}.$$
\end{lemma}

\begin{proof}
Clearly we have
\[
G \cdot Y \subset G \cdot \overline{Y} \subset \overline{G \cdot Y}
\]
and therefore it suffices to show that $G \cdot \overline{Y}$ is
closed in $Z$. The action $\pi\colon G \times Z \to Z$ of $G$ on
$Z$ induces a morphism $\bar\pi\colon G \times^P Z \to Z$ which can
be written as the composition
\[
G \times^P Z \liso G/P \times Z \lto Z.
\]
Here the first morphism is the isomorphism given by $[g,z] \sends
(gP, g\cdot z)$ and the second morphism is the projection. As
$G/P$ is proper, we deduce that the morphism $\bar\pi$ is closed. Now $\overline{Y}$ is
$P$-invariant and therefore $G \times^P \overline{Y}$ is defined, and
it is a closed subscheme of $G \times^P Z$. Therefore $\bar\pi(G
\times^P \overline{Y}) = G\cdot \overline{Y}$ is closed in $Z$.
\end{proof}


The following statements concern images under twisted conjugation:

\begin{theorem}[Lang-Steinberg, see \cite{steinberg}, Theorem 10.1]
\label{thm:lang-steinberg}
Let $G$ be a connected algebraic group and $\phi\colon G\to G$ an isogeny with only a finite number of fixed points. Then the morphism $G\to G$, $g\mapsto g\phi(g)^{-1}$ is surjective.
\end{theorem}

\begin{proposition}\label{lem:Brepr}
Let $G$ be a connected reductive algebraic group with a Borel subgroup $B$ and a maximal torus $T\subset B$.
Let $\phi\colon G\to G$ be an isogeny with $\phi(B)=B$. In (b) also assume that $\phi(T)=T$.
\begin{itemize}
\item[(a)] The morphism $G\times B\to G$, $(g,b)\mapsto gb\phi(g)^{-1}$ is surjective.
\item[(b)] The morphism $G\times T\to G$, $(g,t)\mapsto gt\phi(g)^{-1}$ has dense image.
\end{itemize}
\end{proposition}

If $G$ is simply connected semisimple and $\phi$ is an automorphism of $G$, (b) has been shown by Springer (\cite{Springer:TwistedConjugacy}~Lemma~4).

\begin{proof}
For (a) see \cite{steinberg}, Lemma 7.3. Part (b) and its proof are a slight modification of this.
Equivalently we may show that for some $t_{0}\in T$, the image of the morphism $\tilde\alpha \colon G\times T\to G,  (g,t)\mapsto gtt_{0}\phi(g)^{-1}t_{0}^{-1}$ is dense. For this it will suffice to show that the differential of $\tilde\alpha$ at $1$ is surjective. This differential is the linear map
  \begin{align*}
    \Lie(G)\times \Lie(T) &\to \Lie(G)\\
    (X,Y) &\mapsto X+Y-\Lie(\phi_{t_0})(X),
  \end{align*}
where $\phi_{t_{0}} := \inn(t_{0})\circ \phi$. This linear map has image
  \begin{equation*}
    \Lie(T)+(1-\Lie(\phi_{t_0}))\Lie(G).
  \end{equation*}
Let $B^{-}$ be the Borel subgroup opposite to $B$ with respect to $T$. Since $\phi(B) = B$ and $\phi(T) = T$, the differential of $\phi_{t_{0}}$ at $1$ preserves $\Lie(\ru B)$ and $\Lie(\ru B^{-})$. If we find a $t_{0}$ such that $\Lie(\phi_{t_{0}})$ has no fixed points on $\Lie(\ru B)$ and $\Lie(\ru B^{-})$ we will be done. 

Let $\Phi$ be the set of roots of $G$ with respect to $T$. For each $\alpha\in \Phi$, let $x_{\alpha}$ be a basis vector of $\Lie(U_{\alpha})$, where $U_{\alpha}$ is the unipotent root subgroup of $G$ associated to~$\alpha$. As the isogeny $\phi$ sends $T$ to itself, it induces a bijection $\tilde\phi\colon \Phi\to \Phi$ such that $\phi(U_{\alpha})=U_{\tilde\phi(\alpha)}$. For each $\alpha \in \Phi$ there exists a $c(\alpha)\in k$ such that $\Lie(\phi)(x_{\alpha})=c(\alpha)x_{\tilde\phi(\alpha)}$. This implies $\Lie(\phi_{t_{0}})(x_{\alpha})=\alpha(t_{0})c(\alpha)x_{\tilde\phi(\alpha)}$. Since $\phi_{t_{0}}$ fixes $\ru B$ and $\ru B^{-}$, its differential permutes $\Phi^{+}$ and $\Phi^{-}$, where $\Phi^+$ (resp.~$\Phi^-$) is the set of roots that are positive (resp.~negative) with respect to $B$. Hence $\Lie(\phi_{t_{0}})$ can only have a fixed point in $\Lie(\ru B)$ or $\Lie(\ru B^{-})$ if there exists a cycle $(\alpha_{1},\cdots,\alpha_{n})$ of the permutation $\tilde\phi$ in $\Phi^{+}$ or $\Phi^{-}$ such that
\begin{equation*}
  \prod_{i=1}^{n}\alpha_{i}(t_{0})c(\alpha_{i})=1.
\end{equation*}
This shows that for $t_{0}$ in some non-empty open subset of $T$, the differential $\Lie(\phi_{t_{0}})$ has no fixed points on $\Lie(\ru B)$ and $\Lie(\ru B^{-})$.
\end{proof}


\subsection{Coset Representatives in Coxeter Groups}
\label{Coxeter}

Here we collect some facts about Coxeter groups and root systems which we shall need in the sequel. Let $W$ be a Coxeter group and $S$ its generating set of simple reflections. Let $\ell$ denote the length function on~$W$; thus $\ell(w)$ is the smallest integer $n\ge0$ such that $w=s_1\cdots s_n$ for suitable $s_i\in S$. Any such product with $\ell(w)=n$ is called a reduced expression for~$w$.

Let $I$ be a subset of $S$. We denote by $W_{I}$ the subgroup of $W$ generated by $I$, which is a Coxeter group with set of simple reflections~$I$. Also, we denote by $W^{I}$ (respectively $\leftexp{I}{W}$) the set of elements $w$ of $W$ which have minimal length in their coset $wW_{I}$ (respectively $W_{I}w$). Then every $w \in W$ can be written uniquely as $w=w^{I}\cdot w_{I}=\tilde w_{I}\cdot \leftexp{I}{w}$ with $w_{I}, \tilde w_{I}\in W_{I}$ and $w^{I}\in W^{I}$ and $\leftexp{I}{w}\in\leftexp{I}{W}$. Moreover, these decompositions satisfy $\ell(w)=\ell(w_{I})+\ell(w^{I})=\ell(\tilde w_{I})+\ell(\leftexp{I}{w})$ (see \cite{DDPW}, Proposition 4.16). In particular, $W^{I}$ and $\leftexp{I}{W}$ are systems of representatives for the quotients $W/W_{I}$ and $W_{I}\backslash W$, respectively. The fact that $\ell(w)=\ell(w^{-1})$ for all $w \in W$ implies that $W^{I}=(\leftexp{I}{W})^{-1}$. 

Furthermore, if $J$ is a second subset of $S$, we let $\doubleexp{I}{W}{J}$ denote the set of $x \in W$ which have minimal length in the double coset $W_{I}xW_{J}$. Then $\doubleexp{I}{W}{J}=\leftexp{I}{W}\cap W^{J}$, and it is a system of representatives for $W_{I}\backslash W/W_{J}$ (see \cite{DDPW}~(4.3.2)).

In the next propositions we take an element $x\in \doubleexp{I}{W}{J}$, consider the conjugate subset $\leftexp{x^{-1}}{\!I}\subset W$, and abbreviate $I_x \defeq J\cap\leftexp{x^{-1}}{\!I} \subset J$. Then $\leftexp{I_x}{W_J}$ is the set of elements $w_J$ of $W_J$ which have minimal length in their coset $W_{I_x}w_J$. Likewise $W_I^{I\cap\,\leftexp{x}{\!J}}$ is the set of elements $w_I$ of $W_I$ which have minimal length in their coset $w_IW_{I\cap\,\leftexp{x}{\!J}}$.


\begin{proposition}[Kilmoyer, \cite{DDPW}, Proposition~4.17]
\label{thm:kilmoyer}
  For $x\in \doubleexp{I}{W}{J}$ we have
$$W_{I}\cap\leftexp{x}{W_{J}}=W_{I\cap\,\leftexp{x}{\!J}}
  \qquad\hbox{and}\qquad
  W_J\cap\leftexp{x^{-1}}{W_I}=W_{J\cap\,\leftexp{x^{-1}}{\!I}} = W_{I_x}.$$
\end{proposition}

\begin{proposition}[Howlett, \cite{DDPW}, Proposition 4.18]
\label{prop:howlett}
For any $x \in \doubleexp{I}{W}{J}$, every element $w$ of the double coset $W_{I}xW_{J}$ is uniquely expressible in the form $w=w_{I}xw_{J}$ with $w_{I}\in W_{I}$ and $w_{J}\in \leftexp{I_x}{W_J}$.
Moreover, this decomposition satisfies
  \begin{equation*}
    \ell(w) = \ell(w_{I}xw_{J})=\ell(w_{I})+\ell(x)+\ell(w_{J}).
  \end{equation*}
\end{proposition}

\begin{proposition}\label{howlettcor2}
The set $\leftexp{I}{W}$ is the set of all $xw_J$ for $x\in\doubleexp{I}{W}{J}$ and $w_J\in\leftexp{I_x}{W_J}$.
\end{proposition}

\begin{proof}
Take $x\in\doubleexp{I}{W}{J}$ and $w_J\in\leftexp{I_x}{W_J}$. Then for any $w_I\in \leftexp{I}{W}$, Proposition \ref{prop:howlett} applied to the product $w_Ixw_J$ implies that $\ell(w_Ixw_J)=\ell(w_I)+\ell(x)+\ell(w_J)\geq \ell(x)+\ell(w_J)=\ell(xw_J)$. This proves that $xw_J\in \leftexp{I}{W}$. Conversely take $w\in \leftexp{I}{W}$ and let $w=w_I x w_J$ be its decomposition from Proposition~\ref{prop:howlett}. Then by the first part of the proof we have $xw_J\in \leftexp{I}{W}$. Since $W_Iw=W_Ixw_J$, this implies that $w=xw_J$.
\end{proof}

\begin{proposition}\label{howlettcor2dual}
The set $W^J$ is the set of all $w_Ix$ for $x\in\doubleexp{I}{W}{J}$ and $w_I\in W_I^{I\cap\,\leftexp{x}{\!J}}$.
\end{proposition}

\begin{proof}
Apply Proposition \ref{howlettcor2} with $I$ and $J$ interchanged and invert all elements of~$W^J$.
\end{proof}

Next we recall the Bruhat order on $W$, which we denote by $\leq$ and~$<$. This natural partial order is characterized by the following property: For $x, w\in W$ we have $x\leq w$ if and only if for some (or, equivalently, any) reduced expression $w=s_1\cdots s_n$ as a product of simple reflections $s_i\in S$, one gets a reduced expression for $x$ by removing certain $s_i$ from this product. More information about the Bruhat order can be found in \cite{ccg}, Chapter 2.

Using this order, the set $\leftexp{I}{W}$ can be described as
\UseTheoremCounterForNextEquation
\begin{equation}\label{eq:iwchar}
\leftexp{I}{W}=\{w \in W\mid w < sw \text{ for all } s \in I\}
\end{equation}
(see \cite{ccg}, Definition 2.4.2 and Corollary 2.4.5).

Assume in addition that $W$ is the Weyl group of a root system $\Phi$, with $S$ corresponding to a basis of $\Phi$. Denote the set of positive roots with respect to the given basis by $\Phi^{+}$ and the set of negative roots by $\Phi^{-}$. For $I \subset S$, let $\Phi_{I}$ be the root system spanned by the basis elements corresponding to $I$, and set $\Phi_{I}^{\pm}\defeq\Phi_{I}\cap \Phi^{\pm}$. Then by \cite{carter}, Proposition 2.3.3 we have
\UseTheoremCounterForNextEquation
\begin{equation}\label{eq:wichar1}
W^{I} = \{w\in W\mid w\Phi_{I}^{+}\subset \Phi^{+}\} = \{w\in W\mid w\Phi_{I}^{-}\subset \Phi^{-}\}.
\end{equation}
Also, by \cite{carter}, Proposition 2.2.7, the length of any $w\in W$ is 
\UseTheoremCounterForNextEquation
\begin{equation}\label{eq:lw}
\ell(w) = \#\{\alpha\in\Phi^{+}\mid w\alpha\in\Phi^{-}\}.
\end{equation}

\begin{lemma}\label{lem:refinedlength}
Let $w \in \leftexp{I}{W}$ and write $w = xw_J$ with $x \in \doubleexpIJ{W}$ and $w_J \in W_J$. Then
$$\ell(x) = \#\{\alpha \in \Phi^+ \setminus \Phi_J \mid w\alpha \in \Phi^- \setminus \Phi_I\}.$$
\end{lemma}

\begin{proof}
First note that $\alpha \in \Phi^+$ and $w\alpha \in \Phi^-$ already imply $w\alpha \notin \Phi_I$, because otherwise we would have $\alpha \in w^{-1}\Phi_I^-$, which by \eqref{eq:wichar1} is contained in $\Phi^-$ because $w^{-1} \in W^I$. Thus the right hand side of the claim is $\#\{\alpha \in \Phi^+ \setminus \Phi_J \mid w\alpha \in \Phi^-\}$. Secondly, if $\alpha \in \Phi_J^+$, using again~\eqref{eq:wichar1} and $x \in W^J$ we find that $w\alpha \in \Phi^-$ if and only if $w_J\alpha \in \Phi_J^-$. Thus with~\eqref{eq:lw} we obtain
\begin{eqnarray*}
\#\{\alpha \in \Phi^+ \setminus \Phi_J \mid w\alpha \in \Phi^-\}
&=& \#\{\alpha \in \Phi^+ \mid w\alpha \in \Phi^-\} - \#\{\alpha \in \Phi_J^+ \mid w_J\alpha \in \Phi_J^-\}\\
&=& \ell(w) - \ell(w_J)\ =\ \ell(x).
\end{eqnarray*}
\vskip-15pt\qedhere
\end{proof}


\subsection{Reductive groups, Weyl groups, and parabolics}
\label{Weyl}

Let $G$ be a connected reductive algebraic group, let $B$ be a Borel subgroup of $G$, and let $T$ be a maximal torus of~$B$. 
Let $\Phi(G,T)$ denote the root system of $G$ with respect to~$T$, let $W(G,T) \defeq \Norm_G(T)/T$ denote the associated Weyl group, and let $S(G,B,T) \subset W(G,T)$ denote the set of simple reflections defined by~$B$. Then $W(G,T)$ is a Coxeter group with respect to the subset $S(G,B,T)$. 

A priori this data depends on the pair $(B,T)$. However, any other such pair $(B',T')$ is obtained on conjugating $(B,T)$ by some element $g\in G$ which is unique up to right multiplication by~$T$. Thus conjugation by $g$ induces isomorphisms $\Phi(G,T) \stackrel{\sim}{\to} \Phi(G,T')$ and $W(G,T) \stackrel{\sim}{\to} W(G,T')$ and $S(G,B,T) \stackrel{\sim}{\to} S(G,B',T')$ that are independent of~$g$. Moreover, the isomorphisms associated to any three such pairs are compatible with each other. Thus $\Phi := \Phi(G,T)$ and $W := W(G,T)$ and $S := S(G,B,T)$ for any choice of $(B,T)$ can be viewed as instances of `the' root system and `the' Weyl group and `the' set of simple reflections of $G$, in the sense that up to unique isomorphisms they depend only on~$G$. It then also makes sense to say that the result of a construction (as in Subsection \ref{FirstGw} below) depending on an element of $W$ is independent of $(B,T)$.

For any $w \in W(G,T)$ we fix a representative $\dot w \in \Norm_G(T)$. By choosing representatives attached to a Chevalley system (see \cite{SGA3III}~Exp.~XXIII, \S6) for all $w_1, w_2 \in W$ with $\ell(w_1w_2) = \ell(w_1) + \ell(w_2)$ we obtain
\UseTheoremCounterForNextEquation
\begin{equation}\label{DotMultEq}
\dot w_1 \dot w_2 = (w_1w_2)\dot{\ }.
\end{equation}
In particular the identity element $1\in W$ is represented by the identity element $1\in G$.

A parabolic subgroup of $G$ that contains $B$ is called a \emph{standard parabolic of~$G$}. Any standard parabolic possesses a unique Levi decomposition $P = \ru P \rtimes L$ with $T\subset L$. Any such $L$ is called a standard Levi subgroup of $G$, and the set $I$ of simple reflections in the Weyl group of $L$ is called the \emph{type of~$L$} or \emph{of~$P$}. In this way there is a unique standard parabolic $P_I$ of type $I$ for every subset $I\subset S$, and vice versa. The type of a general parabolic $P$ is by definition the type of the unique standard parabolic conjugate to~$P$; it is independent of $(B,T)$ in the above sense. Any conjugate of a standard Levi subgroup of $G$ is called a Levi subgroup of~$G$.

For any subset $I\subset S$ let $\Par_I$ denote the set of all parabolics of $G$ of type~$I$. Then there is a natural bijection $G/P_I\iso\Par_I$, $gP_I\mapsto\leftexp{g}{P_I}$. For any two subsets $I$, $J\subset S$ we let $G$ act by simultaneous conjugation on $\Par_I\times \Par_J$. As a consequence of the Bruhat decomposition (see \cite{springer} 8.4.6 (3)), the $G$-orbit of any pair $(P,Q)\in \Par_I\times \Par_J$ contains a unique pair of the form $(P_I,\leftexp{\dot x}{P_J})$ with $x\in\doubleexp{I}{W}{J}$. This element $x$ is called the \emph{relative position of $P$ and $Q$} and is denoted by $\relpos(P,Q)$.

We will also use several standard facts about intersections of parabolics and/or Levi subgroups, for instance (see \cite{carter}, Proposition 2.8.9):

\begin{proposition}\label{lem:smallerparabolic}
Let $L$ be a Levi subgroup of~$G$ and $T$ a maximal torus of~$L$. Let $P$ be a parabolic subgroup of $G$ containing~$T$ and $P=\ru P\rtimes H$ its Levi decomposition with $T\subset H$. Then $L \cap P$ is a parabolic subgroup of $L$ with Levi decomposition
$$L \cap P = (L \cap \ru P) \rtimes (L \cap H).$$
If $P$ is a Borel subgroup of $G$, then $L \cap P$ is a Borel subgroup of $L$.
\end{proposition}


\section{Connected algebraic zip data}
\label{AlgZipData}

We now define the central technical notions of this article.

\begin{definition}\label{ZipDatumDef}
A \emph{connected algebraic zip datum} is a tuple $\CZ = (G,P,Q,\varphi)$ consisting of a connected reductive group $G$ with parabolic subgroups $P$ and $Q$ and an isogeny $\varphi\colon P/R_uP\to Q/R_uQ$. The group
\UseTheoremCounterForNextEquation
\begin{equation}\label{ZipGroupDef}
E_\CZ := \bigl\{ (p,q)\in P{\times}Q \bigm| \phi(\pi_{P}(p)) =\pi_Q(q)\bigr\}
\end{equation}
is called the \emph{zip group associated to~$\CZ$}. It acts on $G$ by restriction of the left action
\UseTheoremCounterForNextEquation
\begin{equation}\label{ActionDef}
(P{\times}Q)\times G\to G,\ \ \bigl((p,q),g\bigr)\mapsto pgq^{-1}.
\end{equation}
For any subset $X \subset G$ we denote the union of the $E_\CZ$-orbits of all elements of $X$ by 
\UseTheoremCounterForNextEquation
\begin{equation}\label{OrbDef}
\orb_\CZ(X).
\end{equation}
\end{definition}
Note that if $X$ is a constructible subset of $G$, then so is $\orb_\CZ(X)$. 

Throughout the following sections we fix a connected algebraic zip datum $\CZ = (G,P,Q,\varphi)$. We also abbreviate $U:=\ru P$ and $V:=\ru Q$, so that $\phi$ is an isogeny $P/U\to Q/V$. Our aim is to study the orbit structure of the action of~$E_\CZ$ on $G$. 

\begin{example}\label{PequalGRem}
For dimension reasons we have $P=G$ if and only if $Q=G$. In that case the action of $E_\CZ = \mathop{\rm graph}(\phi)$ is equivalent to the action of $G$ on itself by twisted conjugation $(h,g) \mapsto hg\phi(h)^{-1}$.
\end{example}

In order to work with $\CZ$ it is convenient to fix the following data.

\begin{definition}\label{FrameDef}
A \emph{frame of} $\CZ$ is a tuple $(B,T,g)$ consisting of a Borel subgroup $B$ of~$G$, a maximal torus $T$ of $B$, and an element $g\in G$, such that
\begin{itemize}
\item[(a)] $B\subset Q$,
\item[(b)] $\leftexp{g}{B} \subset P$,
\item[(c)] $\phi\bigl(\pi_P(\leftexp{g}{B})\bigr) = \pi_Q(B)$, and
\item[(d)] $\phi\bigl(\pi_P(\leftexp{g}{T})\bigr) = \pi_Q(T)$.
\end{itemize}
\end{definition}

\begin{proposition}\label{FrameExistence}
Every connected algebraic zip datum possesses a frame.
\end{proposition}

\begin{proof}
Choose a Borel subgroup $B$ of $Q$ and a maximal torus $T$ of~$B$. Let $\bar T'\subset\bar B'\subset P/U$ denote the respective identity components of $\phi^{-1}(\pi_Q(T)) \subset \phi^{-1}(\pi_Q(B)) \subset P/U$. Then $\bar B'$ is a Borel subgroup of $P/U$, and $\bar T'$ is a maximal torus of~$\bar B'$. Thus we have $\bar B' = \pi_P(B')$ for a Borel subgroup $B'$ of~$P$, and $\bar T' = \pi_P(T')$ for some maximal torus $T'$ of~$B'$. Finally take $g \in G$ such that $B' = \leftexp{g}{B}$ and $T' = \leftexp{g}{T}$. Then $(B,T,g)$ is a frame of $\CZ$.
\end{proof}

\begin{proposition}\label{FrameConjugacy}
Let $(B,T,g)$ be a frame of $\CZ$. Then every frame of $\CZ$ has the form $(\leftexp{q}{B},\leftexp{q}{T},pgtq^{-1})$ for $(p,q)\in E_\CZ$ and $t\in T$, and every tuple of this form is a frame of~$\CZ$.
\end{proposition}

\begin{proof}
Let $(B',T',g')$ be another frame of $\CZ$. Since all Borel subgroups of $Q$ are conjugate, we have $B' = \leftexp{q}{B}$ for some element $q\in Q$. Since all maximal tori of $B'$ are conjugate, after multiplying $q$ on the left by an element of $B'$ we may in addition assume that $T'=\leftexp{q}{T}$. Similarly we can find an element $p\in P$ such that $\leftexp{g'}{\!B'}=\leftexp{pg}{B}$ and $\leftexp{g'}{T'}=\leftexp{pg}{T}$. Combining these equations with the defining properties of frames we find that
\begin{eqnarray*}
     \leftexp{\phi(\pi_P(p))}{\pi_Q(B)}
\ =\  \leftexp{\phi(\pi_P(p))}{\phi\bigl(\pi_P(\leftexp{g}{B})\bigr)}
& =& \phi\bigl(\pi_P(\leftexp{pg}{B})\bigr)
\ =\ \phi\bigl(\pi_P(\leftexp{g'}{\!B'})\bigr) \ = \\
&=&  \pi_Q(B')
\ =\ \pi_Q(\leftexp{q}{B})
\ =\ \leftexp{\pi_Q(q)}{\pi_Q(B)},
\end{eqnarray*}
and similarly $\leftexp{\phi(\pi_P(p))}{\pi_Q(T)} = \leftexp{\pi_Q(q)}{\pi_Q(T)}$. Thus $\phi(\pi_P(p)) = \pi_Q(q)\cdot \pi_Q(t')$ for some element $t'\in T$. Since we may still replace $q$ by $qt'$ without changing the above equations, we may without loss of generality assume that $\phi(\pi_P(p)) = \pi_Q(q)$, so that $(p,q) \in E_\CZ$. On the other hand, the above equations imply that $B=\leftexp{g^{-1}p^{-1}g'q}{B}$ and $T=\leftexp{g^{-1}p^{-1}g'q}{T}$, so that $t := g^{-1}p^{-1}g'q \in T$ and hence $g'=pgtq^{-1}$. This proves the first assertion. The second involves a straightforward calculation that is left to the conscientious reader.
\end{proof}

Throughout the following sections we fix a frame $(B,T,g)$ of~$\CZ$. It determines unique Levi components $\leftexp{g}{T} \subset L\subset P$ and $T\subset M\subset Q$. Via the isomorphisms $L\stackrel{\sim}{\to}P/U$ and $M\stackrel{\sim}{\to}Q/V$ we can then identify $\phi$ with an isogeny $\phi\colon L\to M$. The zip group then becomes
\UseTheoremCounterForNextEquation
\begin{equation}\label{ZipDecomp}
E_\CZ\ =\ \bigl\{(u\ell,v\phi(\ell)) \,\bigm|\, u\in U,\, v\in V,\, \ell\in L \bigr\}
\end{equation}
and acts on $G$ by $((u\ell,v\phi(\ell)) , g) \sends u\ell g\phi(\ell)^{-1}v^{-1}$.
Moreover, conditions \ref{FrameDef} (c) and (d) are then equivalent to
\UseTheoremCounterForNextEquation
\begin{equation}\label{LeviFrameCond}
\varphi(\leftexp{g}{B} \cap L) = B \cap M, \qquad\hbox{and}\qquad \phi(\leftexp{g}{T}) = T,
\end{equation}
which are a Borel subgroup and a maximal torus of~$M$, respectively.

Let $\Phi$ be the root system, $W$ the Weyl group, and $S\subset W$ the set of simple reflections of $G$ with respect to $(B,T)$. Let $I \subset S$ be the type of $\leftexp{g^{-1}}{\!P}$ and $J \subset S$ the type of~$Q$. Then $M$ has root system $\Phi_J$, Weyl group~$W_J$, and set of simple reflections $J\subset W_J$. Similarly $\leftexp{g^{-1}}{\!L}$ has root system $\Phi_I$, Weyl group~$W_I$, and set of simple reflections $I\subset W_I$, and the inner automorphism $\inn(g)$ identifies these with the corresponding objects associated to~$L$. Moreover, the equations \eqref{LeviFrameCond} imply that $\phi\circ\inn(g)$ induces an isomorphism of Coxeter systems
\begin{equation}\label{DefinePsi}
\psi\colon (W_I,I) \iso (W_J,J).
\end{equation}

Recall that $\Phi$, $W$, and $S$ can be viewed as independent of the chosen frame, as explained in Subsection~\ref{Weyl}.

\begin{proposition}\label{Indep1}
The subsets $I$, $J$ and the isomorphism $\psi$ are independent of the frame.
\end{proposition}

\begin{proof}
Consider another frame $(\leftexp{q}{B},\leftexp{q}{T},pgtq^{-1})$ with $(p,q)\in E_\CZ$ and $t\in T$, as in Proposition~\ref{FrameConjugacy}. Then we have a commutative diagram
$$\xymatrix@C+40pt{
(\leftexp{g^{-1}}{\!L},B,T) \ar[d]_{\inn(qt^{-1})}^\wr \ar[r]^-{\inn(g)}_\sim &
(L,\leftexp{g}{B},\leftexp{g}{T}) \ar[d]_{\inn(p)}^\wr \ar[r]^-{\phi} &
(M,B,T) \ar[d]_{\inn(q)}^\wr \\
(\leftexp{qg^{-1}}{\!L},\leftexp{q}{B},\leftexp{q}{T}) \ar[r]^-{\inn(pgtq^{-1})}_\sim &
(\leftexp{p}{L},\leftexp{pg}{B},\leftexp{pg}{T}) \ar[r]^-{\phi} &
(\leftexp{q}{M},\leftexp{q}{B},\leftexp{q}{T}) \rlap{,} \\}$$
whose upper row contains the data inducing $\psi$ for the old frame and whose lower row is the analogue for the new frame. Since the vertical arrows are inner automorphisms, they induce the identity on the abstract Coxeter system $(W,S)$ of $G$ as explained in Subsection~\ref{Weyl}. Everything follows from this.
\end{proof}


\section{Induction step}\label{Bruhat}

We keep the notations of the preceding section. Since $\leftexp{g^{-1}}{\!P}$ and $Q$ are parabolic subgroups containing the same Borel $B$, by Bruhat (see \cite{springer} 8.4.6 (3)) we have a disjoint decomposition
$$G\ =\ \coprod_{x\in\doubleexpIJ{W}} \leftexp{g^{-1}}{\!P} \dot x Q.$$
Left translation by $g$ turns this into a disjoint decomposition
\UseTheoremCounterForNextEquation
\begin{equation}\label{BruhatDecomp}
G\ =\ \coprod_{x\in\doubleexpIJ{W}} P g \dot x Q.
\end{equation}
Here each component $Pg\dot xQ$ is an irreducible locally closed subvariety of $G$ that is invariant under the action of~$E_\CZ$. In this section we fix an element $x\in\doubleexpIJ{W}$ and establish a bijection between the $E_\CZ$-orbits in $Pg\dot xQ$ and the orbits of another zip datum constructed from $\CZ$ and~$\dot x$. This will allow us to prove facts about the orbits inductively. The base case of the induction occurs when the decomposition possesses just one piece, i.e., when $P=Q=G$.

\begin{lemma}\label{step1}
The stabilizer of $g\dot xQ \subset Pg\dot xQ$ in $E_\CZ$ is the subgroup
$$E_{\CZ,\dot x} := \bigl\{ (p,q)\in E_\CZ \bigm| p\in P\cap \leftexp{g\dot x}{Q} \bigr\},$$
and the action of $E_\CZ$ induces an $E_\CZ$-equivariant isomorphism
$$E_\CZ\times^{E_{\CZ,\dot x}} g\dot xQ \stackrel{\sim}{\longto} Pg\dot xQ,\ 
[((p,q),h)] \mapsto phq^{-1}.$$
\end{lemma}

\begin{proof}
The action \eqref{ActionDef} of $(p,q)\in E_\CZ$ on $Pg\dot xQ$ induces the action on the quotient $Pg\dot xQ/Q$ by left multiplication with~$p$. From \eqref{ZipGroupDef} we see that the latter action is transitive, and the stabilizer of the point $g\dot xQ$ is $E_{\CZ,\dot x}$; hence there is an $E_\CZ$-equivariant isomorphism $Pg\dot xQ/Q \cong E_\CZ / E_{\CZ,\dot x}$. Thus everything follows by applying Lemma~\ref{lem:pushout} to the projection morphism $Pg\dot xQ \onto Pg\dot xQ/Q \cong E_\CZ / E_{\CZ,\dot x}$.
\end{proof}

\begin{construction}
\label{konstr:zn}
Consider the following subgroups of the connected reductive algebraic group~$M$ (which are independent of the representative $\dot x$ of~$x$):
$$\begin{tabular}{ll}
\llap{$P_x$}$\;\defeq M\cap \leftexp{\dot x^{-1}g^{-1}}{\!P}$,\qquad\quad &
\llap{$Q_x$}$\;\defeq\phi(L\cap\leftexp{g\dot x}{Q})$, \\
\llap{$U_x$}$\;\defeq M\cap \leftexp{\dot x^{-1}g^{-1}}{\!U}$, &
\llap{$V_x$}$\;\defeq\phi(L\cap\leftexp{g\dot x}{V})$, \\
\llap{$L_x$}$\;\defeq M\cap \leftexp{\dot x^{-1}g^{-1}}{\!L}$, &
\llap{$M_x$}$\;\defeq\phi(L\cap\leftexp{g\dot x}{M})$. \\
\end{tabular}$$
Proposition~\ref{lem:smallerparabolic} shows that $P_x$ is a parabolic with unipotent radical $U_x$ and Levi component~$L_x$, and that $Q_x$ is a parabolic with unipotent radical $V_x$ and Levi component~$M_x$. Moreover, $\phi\circ\inn(g\dot x)$ induces an isogeny $\phi_{\dot x}\colon L_x\to M_x$, or equivalently $P_x/U_x\to Q_x/V_x$. Thus we obtain a connected algebraic zip datum $\CZ_{\dot x} \defeq (M,P_x,Q_x,\phi_{\dot x})$. By \eqref{ZipDecomp} its zip group is
\UseTheoremCounterForNextEquation
\begin{equation}\label{ZipGroupW}
E_{\CZ_{\dot x}}\ =\ \bigl\{(u'\ell',v'\phi_{\dot x}(\ell')) 
                  \,\bigm|\, u'\in U_x,\, v'\in V_x,\, \ell'\in L_x \bigr\}.
\end{equation}
\end{construction}

\begin{lemma}\label{step2}
There is a surjective homomorphism 
$$E_{\CZ,{\dot x}} \onto E_{\CZ_{\dot x}},\ (p,q)\mapsto(m,\phi(\ell)),$$
where $p=u\ell$ for $u\in U$ and $\ell\in L$, and $\leftexp{\dot x^{-1}g^{-1}}{\!p}=vm$ for $v\in V$ and $m\in M$.
\end{lemma}

\begin{proof}
For ease of notation abbreviate $h:=g\dot x$, so that $\leftexp{h}{T} = \leftexp{g}{T} \subset L$ and therefore $T\subset\leftexp{h^{-1}}{L} \subset\leftexp{h^{-1}}{P}$. Thus $\leftexp{h^{-1}}{\!P}$ and $Q$ are parabolics of $G$ with $T$-invariant Levi decompositions $\leftexp{h^{-1}}{\!P} = \leftexp{h^{-1}}{\!U} \rtimes \leftexp{h^{-1}}{\!L}$ and $Q=V\rtimes M$. It follows (see \cite{carter} Thm. 2.8.7) that any element of $\leftexp{h^{-1}}{\!P} \cap Q$ can be written as a product $abu'\ell'$ with unique
$$\begin{tabular}{ll}
\llap{$a$}$\;\in \leftexp{h^{-1}}{\!U} \cap V$,\qquad\quad &
\llap{$u'$}$\;\in\leftexp{h^{-1}}{\!U} \cap M = U_x$, \\
\llap{$b$}$\;\in \leftexp{h^{-1}}{\!L} \cap V$, &
\llap{$\ell'$}$\;\in\leftexp{h^{-1}}{\!L} \cap M = L_x$. \\
\end{tabular}$$
Consider $(p,q) \in E_{\CZ,{\dot x}}$ with $p=u\ell$ and $\leftexp{h^{-1}}{\!p}=vm$ as in the lemma. Then we can write the element $\leftexp{h^{-1}}{\!p} = abu'\ell' \in \leftexp{h^{-1}}{\!P} \cap Q$ in the indicated fashion. Comparing the different factorizations yields the equations $v=ab$, $m=u'\ell'$, $u=\leftexp{h}{(abu'}b^{-1})$, and $\ell=\leftexp{h}{(b\ell')}$. Thus $\phi(\ell) = \phi(\leftexp{h}{b})\phi(\leftexp{h}{\ell'}) = v'\phi_{\dot x}(\ell')$ with $v' := \phi(\leftexp{h}{b}) \in \phi(L\cap\leftexp{h}{V}) = V_x$. In view of \eqref{ZipGroupW} it follows that $(m,\phi(\ell)) = (u'\ell',v'\phi_{\dot x}(\ell'))$ lies in $E_{\CZ_{\dot x}}$, and so the map in question is well-defined. Since $m$ and $\ell$ are obtained from $p$ by projection to Levi components, the map is a homomorphism. Conversely, every element of $E_{\CZ_{\dot x}}$ can be obtained in this way from some element $p\in P\cap\leftexp{h}{Q}$. By \eqref{ZipDecomp} we can then also find $q\in Q$ with $(p,q) \in E_{\CZ,{\dot x}}$. Thus the map is surjective, and we are done.
\end{proof}

\begin{lemma}\label{step3}
The surjective morphism 
$$\pi\colon g\dot xQ \onto M,\ g\dot x\tilde m\tilde v \mapsto\tilde m$$
for $\tilde m\in M$ and $\tilde v\in V$ is equivariant under the group $E_{\CZ,{\dot x}}$, which acts on $g\dot xQ$ as in Lemma \ref{step1} and on $M$ through the homomorphism from Lemma~\ref{step2}.
\end{lemma}

\begin{proof}
Take $(p,q)\in E_{\CZ,{\dot x}}$ with $p=u\ell$ and $\leftexp{\dot x^{-1}g^{-1}}{\!p}=vm$ as in Lemma~\ref{step2}. Then \eqref{ZipDecomp} implies that $q=v_1\phi(\ell)$ for some $v_1\in V$. Thus the action of $(p,q)$ sends $g\dot x\tilde m\tilde v \in g\dot xQ$ to the element
$$pg\dot x\cdot\tilde m\tilde v\cdot q^{-1} = g\dot x vm\cdot\tilde m\tilde v\cdot\phi(\ell)^{-1}v_1^{-1}
 = g\dot x \cdot m\tilde m \phi(\ell)^{-1} \cdot \bigl(\hbox{an element of $V$}\bigr).$$
The morphism $\pi$ maps this element to $m\tilde m \phi(\ell)^{-1} \in M$. But this is also the image 
of $\tilde m = \pi(g\dot x\tilde m\tilde v)$ under the action of $(m,\phi(\ell)) \in E_{\CZ_{\dot x}}$. Thus the morphism is equivariant.
\end{proof}

\begin{proposition}\label{reduction}
There is a closure-preserving bijection between $E_{\CZ_{\dot x}}$-invariant subsets $Y\subset M$ and $E_{\CZ}$-invariant subsets $X\subset Pg\dot xQ$, defined by $Y = M\cap \dot x^{-1}g^{-1}X$ and $X=\orb_\CZ(g\dot xY)$. Moreover, $Y$ is a subvariety if and only if $X$ is one, and in that case $X\cong E_\CZ\times^{E_{\CZ,{\dot x}}} \pi^{-1}(Y)$.
\end{proposition}

\begin{proof}
From \eqref{ZipGroupDef} and \eqref{ActionDef} we see that the subgroup $V \cong \{(1,v)\mid v\in V\} \subset E_{\CZ,{\dot x}}$ acts by right translation on $g\dot xQ$. Thus every $E_{\CZ,{\dot x}}$-invariant subset of $g\dot xQ$ is a union of cosets of $V$ and therefore of the form $Z=g\dot xYV=\pi^{-1}(Y)$ for a subset $Y\subset M$, which moreover satisfies $Y = M\cap \dot x^{-1}g^{-1}Z$. By Lemmas \ref{step2} and \ref{step3} this defines a bijection between $E_{\CZ_{\dot x}}$-invariant subsets $Y\subset M$ and $E_{\CZ,{\dot x}}$-invariant subsets $Z\subset g\dot xQ$. On the other hand, Lemmas \ref{lem:pushoutorbits} and \ref{step1} yield a bijection between $E_{\CZ,{\dot x}}$-invariant subsets $Z\subset g\dot xQ$ and $E_{\CZ}$-invariant subsets $X\subset Pg\dot xQ$ that is characterized by $Z=g\dot xQ\cap X$ and $X=\orb_\CZ(Z)$. Together we obtain the desired bijection with $Y = M\cap \dot x^{-1}g^{-1}(g\dot xQ\cap X) = M\cap \dot x^{-1}g^{-1}X$ and $X=\orb_\CZ(g\dot xYV)=\orb_\CZ(g\dot xY)$.

The equations $Z=\pi^{-1}(Y)$ and $Y = M\cap \dot x^{-1}g^{-1}Z$ imply that the bijection between $Y$ and $Z$ preserves closures and maps subvarieties to subvarieties. The corresponding facts for the bijection between $Z$ and $X$ follow from Lemma~\ref{lem:pushoutorbits}, which also implies the last statement.
\end{proof}

\begin{proposition}\label{reductiondimension1}
If $X$ and $Y$ in Proposition \ref{reduction} are subvarieties, then
$$\dim X\ =\ \dim Y + \dim P - \dim P_x + \ell(x).$$
\end{proposition}

\begin{proof}
By the definition of $E_{\CZ,{\dot x}}$ we have
$$\dim E_\CZ - \dim E_{\CZ,{\dot x}} \ =\ \dim P - \dim(P\cap\leftexp{g\dot x}{Q})
\ =\ \dim P - \dim P_x - \dim(P\cap\leftexp{g\dot x}{V}).$$
With the last statement of Proposition~\ref{reduction} this implies that 
$$\dim X 
\ =\ \dim Y + \dim V + \dim P - \dim P_x - \dim(P\cap\leftexp{g\dot x}{V}).$$
From the decomposition of $V$ into root subgroups it follows that $\dim V - \dim(P\cap\leftexp{g\dot x}{V}) = \dim V - \dim(V\cap\leftexp{\dot x^{-1}g^{-1}}{\!P})$ is the cardinality of the set
$$\{\alpha\in\Phi^{+}\setminus \Phi_{J}\mid x\alpha\in \Phi^{-}\setminus\Phi_{I}\}.$$
By Lemma \ref{lem:refinedlength} for $w_J=1$ this cardinality is $\ell(x)$.
\end{proof}

\begin{lemma}\label{ozozw}
For any subset $Y\subset M$ we have $\orb_\CZ(g\dot x\orb_{\CZ_{\dot x}}(Y)) = \orb_\CZ(g\dot xY)$.
\end{lemma}

\begin{proof}
It suffices to show that $g\dot x\orb_{\CZ_{\dot x}}(Y) \subset \orb_\CZ(g\dot xY)$, which follows from a straightforward calculation that is left to the reader. Alternatively the formula can be deduced from the formal properties stated in Proposition~\ref{reduction}.
\end{proof}

We can also give an inductive description of the stabilizers of points in $Pg\dot xQ$. However, this does not give the scheme-theoretic stabilizers, which may in fact be non-reduced. Likewise, the following lemma does not describe the scheme-theoretic kernel:

\begin{lemma}\label{step4}
The kernel of the homomorphism from Lemma \ref{step2} is $(U\cap\leftexp{g\dot x}{V})\times V$.
\end{lemma}

\begin{proof}
Let $p=u\ell$ and $\leftexp{\dot x^{-1}g^{-1}}{\!p}=vm$ be as in Lemma \ref{step2}. Then $(p,q)$ is in the kernel if and only if $m=1$ and $\phi(\ell)=1$. The first equation is equivalent to $p=\leftexp{g\dot x}{v} \in \leftexp{g\dot x}{V}$, which implies that $\ell$ is unipotent. Being in the kernel of the isogeny $\phi$ is then equivalent to $\ell=1$. Thus the second equation is equivalent to $p\in U$, and the two together are equivalent to $p\in U\cap\leftexp{g\dot x}{V}$. By \eqref{ZipDecomp} we then have $q\in V$, and so we are done.
\end{proof}

\begin{proposition}\label{stabilizerreduction}
For any $m\in M$ there is a short exact sequence
$$1\longto U\cap\leftexp{g\dot x}{V}\longto \Stab_{E_\CZ}(g\dot xm) 
\stackrel{\ref{step2}}{\longto} \Stab_{E_{\CZ_{\dot x}}}(m)\longto 1.$$
\end{proposition}

\begin{proof}
The second half of Lemma \ref{step1} and Lemma \ref{step3} imply that we have an equality, respectively a homomorphism
$$\Stab_{E_\CZ}(g\dot xm)
\ =\ \Stab_{E_{\CZ,{\dot x}}}(g\dot xm)
\stackrel{\ref{step2}}{\longto} \Stab_{E_{\CZ_{\dot x}}}(m).$$
This homomorphism is surjective, because the subgroup $V \cong \{(1,v)\mid v\in V\}$ contained in the kernel of the surjection $E_{\CZ,{\dot x}} \onto E_{\CZ_{\dot x}}$ acts transitively on the fibers of~$\pi$. By Lemma \ref{step4} the kernel is the stabilizer of $g\dot xm$ in the group $(U\cap\leftexp{g\dot x}{V})\times V$ acting by left and right translation. This stabilizer consists of $(u,\leftexp{(g\dot xm)^{-1}}{\!u})$ for all $u\in U\cap\leftexp{g\dot x}{V}$, and we are done.
\end{proof}


Finally, the assumption $x\in\doubleexpIJ{W}$ allows us to construct a frame of~$\CZ_{\dot x}$:

\begin{proposition}\label{ReductionFrame}
The tuple $(M\cap B,T,1)$ is a frame of $\CZ_{\dot x}$, and the associated Levi components of $P_x$ and $Q_x$ are $L_x$ and $M_x$, respectively.
\end{proposition}

\begin{proof}
First, the assumptions $T\subset M$ and $\leftexp{g\dot x}{T} = \leftexp{g}{T} \subset L$ imply that $T\subset M\cap\leftexp{\dot x^{-1}g^{-1}}{\!L}$, the latter being $L_x$ by Construction~\ref{konstr:zn}. Together with the equation $\phi(\leftexp{g}{T}) = T$ from \eqref{LeviFrameCond} they also imply that $T = \phi(\leftexp{g\dot x}{T}) \subset \phi(L\cap\leftexp{g\dot x}{M})$, the latter being $M_x$ by Construction~\ref{konstr:zn}. This proves the statement about the Levi components. We can also directly deduce that $\phi_{\dot x}(T) = \phi(\leftexp{g\dot x}{T}) = T$.

Next, as $T$ is a common maximal torus of $M$ and $B$, Proposition \ref{lem:smallerparabolic} implies that $M\cap B$ is a Borel subgroup of~$M$. Recall that $M$ has the root system~$\Phi_J$, so that $M\cap B$ corresponds to the subset $\Phi_J^+ = \Phi_J\cap\Phi^+$. For the same reasons $M\cap \leftexp{\dot x^{-1}}{\!B}$ is a Borel subgroup of $M$ corresponding to the subset $\Phi_J\cap x^{-1}\Phi^+$. But with \eqref{eq:wichar1} the assumption $x\in\doubleexpIJ{W} \subset W^J$ implies that $x\Phi_J^+\subset \Phi^+$, and hence $\Phi_J^+ \subset \Phi_J\cap x^{-1}\Phi^+$. Since both subsets correspond to Borel subgroups, they must then coincide, and therefore $M\cap B = M\cap \leftexp{\dot x^{-1}}{\!B}$. With the inclusion $\leftexp{g}{B}\subset P$ from \eqref{FrameDef} we deduce that
$$M\cap B \ =\ 
  M\cap \leftexp{\dot x^{-1}}{\!B} \ \subset\ 
  M\cap \leftexp{\dot x^{-1}g^{-1}}{\!P}
  \;\stackrel{\ref{konstr:zn}}{=}\; P_x.$$
In the same way one shows that $L\cap\leftexp{g}{B} = L\cap \leftexp{ g\dot x}{B}$, which together with $B\subset Q$ implies that
$$M\cap B \,\stackrel{\eqref{LeviFrameCond}}{=}\,
  \phi(L\cap\leftexp{g}{B}) \ =\ 
  \phi(L\cap \leftexp{ g\dot x}{B}) \ \subset\ 
  \phi(L\cap \leftexp{ g\dot x}{Q}) 
  \;\stackrel{\ref{konstr:zn}}{=}\; Q_x.$$
The equation $M\cap B = \phi(L\cap \leftexp{ g\dot x}{B})$ and Construction \ref{konstr:zn} also imply that
$$\phi_{\dot x}\bigl((M\cap B)\cap L_x\bigr) \ =\ 
\phi\bigl(\leftexp{ g\dot x}{M}\cap\leftexp{ g\dot x}{B}\cap L\bigr) \ \subset\ 
\phi\bigl(L\cap\leftexp{ g\dot x}{M}\bigr) \cap\phi\bigl(L\cap\leftexp{ g\dot x}{B}\bigr) \ =\ 
(M\cap B) \cap M_x.$$
As both sides of this inclusion are Borel subgroups of $M_x$, they must be equal. Thus $(M\cap B,T,1)$ satisfies Definition \ref{FrameDef} in the variant \eqref{LeviFrameCond}, as desired.
\end{proof}

Recall that $M$ has the Weyl group $W_J$ with the set of simple reflections~$J$, and that $\psi\colon {W_I \iso W_J}$ is the isomorphism induced by $\phi\circ\inn(g)$.

\begin{proposition}\label{nicereduction1}
\begin{itemize}
\item[(a)] The type of the parabolic $P_x$ of $M$ is $I_x\defeq J\cap\leftexp{x^{-1}}{\!I}$.
\item[(b)] The type of the parabolic $Q_x$ of $M$ is $J_x\defeq \psi(I\cap\leftexp{x}{\!J})$.
\item[(c)] The isomomorphism $\psi_x\colon W_{I_x}\iso W_{J_x}$ induced by $\phi_{\dot x}$ is the restriction of $\psi\circ\inn(x)$.
\end{itemize}
\end{proposition}

\begin{proof}
Proposition \ref{thm:kilmoyer} implies that $L_x=M\cap\leftexp{\dot x^{-1}g^{-1}}{\!L}$ has the Weyl group $W_J\cap\leftexp{x^{-1}}{W_I}=W_{I_x}$, which shows~(a). Likewise $M_x=\phi(L\cap\leftexp{g\dot x}{M})$ has the Weyl group $\psi(W_{I}\cap\leftexp{x}{W_{J}})=W_{J_x}$, which implies~(b). Finally, (c) follows from $\phi_{\dot x} = \phi\circ\inn(g\dot x)$.
\end{proof}


\section{Decomposition of $G$}
\label{Gw}

In this section we construct a natural decomposition of $G$ into finitely many $E_\CZ$-invariant subvarieties~$G^w$.


\subsection{The Levi subgroup $H_w$}
\label{HwSec}

Fix an element $w \in \leftexp{I}{W}$. Note that we can compare any subgroup $H$ of $\leftexp{\dot w^{-1}g^{-1}}{\!L}$ with its image $\phi\circ\inn(g\dot w)(H)$ in~$M$, because both are subgroups of~$G$. Moreover, the collection of all such $H$ satisfying $\phi\circ\inn(g\dot w)(H) = H$ possesses a unique largest element, namely the subgroup generated by all such subgroups.

\begin{definition}
\label{def:HwMult} 
We let $H_w$ denote the unique largest subgroup of $\leftexp{\dot w^{-1}g^{-1}}{\!L}$ satisfying $\phi\circ\inn(g\dot w)(H_w) = H_w$. We let $\phi_{\dot w}\colon H_w\to H_w$ denote the isogeny induced by ${\phi\circ\inn(g\dot w)}$, and let $H_w$ act on itself from the left by the twisted conjugation $(h,h') \mapsto hh'\phi_{\dot w}(h)^{-1}$. 
\end{definition}

\begin{remark}\label{kramer}
Since $\phi\circ\inn(g\dot w)(T) = \phi(\leftexp{g}{T}) = T$ by~\eqref{LeviFrameCond}, the defining property of $H_w$ implies that $T\subset H_w$. Thus $H_w$ does not depend on the choice of representative $\dot w$ of~$w$, justifying the notation~$H_w$. Also, in the case that $w=x\in\doubleexpIJ{W}$ observe that the $\phi_{\dot w}$ defined here is the restriction to $H_w$ of the isogeny $\phi_{\dot x}$ from Construction~\ref{konstr:zn}. Using the same notation for both is therefore only mildly abusive.
\end{remark}

\begin{example}\label{PequalGRem3}
In the case $P=Q=G$ from Example~\ref{PequalGRem} we have $M=L=G$ and $I=J=\psi(J)=S$ and hence $\leftexp{I}{W} = \{1\}$ and $H_1=G$.
\end{example}

To analyze $H_w$ in the general case we apply the induction step from Section~\ref{Bruhat}. Let $w=xw_J$ be the decomposition from Proposition \ref{howlettcor2} with $x\in\doubleexpIJ{W}$ and $w_J\in \leftexp{I_x}{W_J}$ for $I_x = J\cap\leftexp{x^{-1}}{\!I}$. Since $W_J$ is the Weyl group of $M$, and $I_x$ is the type of the parabolic $P_x\subset M$ by Proposition \ref{nicereduction1}~(a), we can also apply Definition \ref{def:HwMult} to the pair $(\CZ_{\dot x},w_J)$ in place of $(\CZ,w)$.

\begin{lemma}
\label{lem:iwreduction2}
The subgroup $H_w$ and the isogeny $\phi_{\dot w}$ associated to $(\CZ,w)$ in Definition \ref{def:HwMult} are equal to those associated to $(\CZ_{\dot x},w_J)$.
\end{lemma}

\begin{proof}
Since $\dot w_J\in M=\phi(L)$, Definition \ref{def:HwMult} and Construction \ref{konstr:zn} imply that
$$H_w\ \subset\ M \cap \leftexp{\dot w_J^{-1}\dot x^{-1}g^{-1}}{\!L}
\ =\ \leftexp{\dot w_J^{-1}}{\!\!\bigl(}M\cap\leftexp{\dot x^{-1}g^{-1}}{\!L}\bigr)
\ =\ \leftexp{\dot w_J^{-1}}{\!L_x}$$
and that $\phi_{\dot x}\circ\inn(\dot w_J)(H_w) = \phi\circ\inn(g\dot w)(H_w) = H_w$. Since $H_w$ is the largest subgroup of $\leftexp{\dot w^{-1}g^{-1}}{\!L}$ with this property, it is also the largest in $\leftexp{\dot w_J^{-1}}{\!L_x}$.
\end{proof}

\begin{remark}
The preceding lemma implies that $H_w$ and $\phi_{\dot w}$ also remain the same if we repeat the induction step with $(\CZ_{\dot x},w_J)$ in place of $(\CZ,w)$, and so on. When the process becomes stationary, we have reached a pair consisting of a zip datum as in Example~\ref{PequalGRem3} and the Weyl group element~$1$, whose underlying connected reductive group and isogeny are $H_w$ and~$\phi_{\dot w}$. This induction process is the idea underlying many proofs throughout this section.
\end{remark}

\begin{proposition}\label{Kw5}
The subgroup $H_w$ is the standard Levi subgroup of $G$ containing $T$ whose set of simple reflections is the unique largest subset $K_w$ of $\leftexp{w^{-1}}{\!I}$ satisfying $\psi\circ\inn(w)(K_w) = K_w$.
\end{proposition}

\begin{proof}
For any subset $K$ of $\leftexp{w^{-1}}{\!I}$ the equality $\psi\circ\inn(w)(K)=K$ makes sense, because both sides are subsets of~$W$. The collection of all such $K$ satisfying that equality possesses a unique largest element~$K_w$, namely the union of all of them. Then $K_w = \psi\circ\inn(w)(K_w) \subset \psi(I)=J\subset S$, and so $K_w$ consists of simple reflections. 

Let $H$ denote the standard Levi subgroup of $G$ containing $T$ with the set of simple reflections~$K_w$. Then the isogeny $\phi\circ\inn(g\dot w)\colon \leftexp{\dot w^{-1}g^{-1}}{\!L} \to M$ sends $T$ to itself by Remark \ref{kramer}, and the associated isomorphism of Weyl groups $\psi\circ\inn(w)\colon {\leftexp{w^{-1}}{W_I} \to W_J}$ sends $K_w$ to itself by construction. Together this implies that $\phi\circ\inn(g\dot w)(H) = H$ and hence $H\subset H_w$.

We now prove the equality $H_w=H$ by induction on $\dim G$. In the base case $M=G$ we have  $I=J=S$ and $w=1$ and thus $K_1=S$ and $H=G$, while $H_1=G$ by Example~\ref{PequalGRem3}; hence we are done. Otherwise write $w=xw_J$ as above. Then Lemma~\ref{lem:iwreduction2} and the induction hypothesis show that $H_w$ is a Levi subgroup of $M$ containing $T$ with a set of simple reflections $K\subset \leftexp{w_J^{-1}}{\!I_x}$ satisfying $\psi_x\circ\inn(w_J)(K) = K$. But $\leftexp{w_J^{-1}}{\!I_x} = \leftexp{w_J^{-1}}{\!(}J\cap\leftexp{x^{-1}}{\!I}) \subset \leftexp{w^{-1}}{\!I}$ and $\psi_x\circ\inn(w_J)$ is the restriction of $\psi\circ\inn(x)\circ\inn(w_J)=\psi\circ\inn(w)$. By the maximality of~$K_w$ we thus have $K\subset K_w$ and therefore $H_w\subset H$. Together with the earlier inequality $H\subset H_w$ we deduce that $H_w=H$, as desired.
\end{proof}


\subsection{First description of $G^w$}
\label{FirstGw}

\begin{definition}\label{GwDef1}
For any $w \in \leftexp{I}{W}$ we set $G^{w} \defeq \orb_\CZ(g\dot wH_w)$.
\end{definition}

\begin{proposition}\label{GwIndependentProp}
The set $G^w$ does not depend on the representative $\dot w$ of $w$ or the frame.
\end{proposition}

\begin{proof}
The independence of $\dot w$ follows from the inclusion $T\subset H_w$. For the rest note first that by Propositions \ref{Indep1} and \ref{Kw5} the set $K_w$ is independent of the frame. Consider another frame $(\leftexp{q}{B},\leftexp{q}{T},pgtq^{-1})$ for $(p,q)\in E_\CZ$ and $t\in T$, as in Proposition~\ref{FrameConjugacy}. Recall from Subsection \ref{Weyl} that the isomorphism $W(G,T) \iso W(G,\leftexp{q}{T})$ is induced by $\inn(q): \Norm_G(T) \iso \Norm_G(\leftexp{q}{T})$. It follows that $w\in\leftexp{I}{W}$ as an element of the abstract Weyl group of $G$ is represented by $q\dot wq^{-1}\in \Norm_G(\leftexp{q}{T})$, and with Proposition \ref{Kw5} it follows that the Levi subgroup associated to $w$ and the new frame is~$\leftexp{q}{H_w}$. Thus the right hand side in Definition \ref{GwDef1} associated to the new frame is
$$\orb_\CZ\bigl((pgtq^{-1})(q\dot wq^{-1})\,\leftexp{q}{\!H_w}\bigr)
= \orb_\CZ\bigl(p g t \dot w H_wq^{-1}\bigr)
= \orb_\CZ\bigl(g t \dot w H_w\bigr)
= \orb_\CZ\bigl(g \dot w H_w\bigr),$$
where the second equation follows from $(p,q)\in E_\CZ$ and the third from $\dot w^{-1}t\dot w\in T\subset H_w$. Thus $G^w$ is independent of the frame.
\end{proof}

In Example \ref{PequalGRem3} we have $H_1=G$ and hence $G^1=G$. Otherwise recall from Proposition \ref{ReductionFrame} that $\CZ_{\dot x}$ has the frame $(M\cap B,T,1)$. Thus by Lemma \ref{lem:iwreduction2}, the subset associated to $(\CZ_{\dot x},w_J)$ by Definition \ref{GwDef1} is $M^{w_J} \defeq \orb_{\CZ_{\dot x}}(\dot w_JH_w)$.

\begin{lemma}\label{lem:nicereduction2}
Under the bijection of Proposition \ref{reduction}, the subset $M^{w_J}\subset M$ corresponds to the subset $G^w\subset Pg\dot xQ$. In particular $G^w = \orb_\CZ(g\dot xM^{w_J})$. Also, there is a bijection between the $E_{\CZ_{\dot x}}$-orbits $X'\subset M^{w_J}$ and the $E_\CZ$-orbits $X\subset G^w$, defined by $X=\orb_\CZ(g\dot xX')$.
\end{lemma}

\begin{proof}
Using, in this order, the definition of $G^w$, the equation \eqref{DotMultEq}, Lemma \ref{ozozw}, and the definition of $M^{w_J}$ we find that
$$G^w = \orb_\CZ(g\dot w H_w) = \orb_\CZ(g\dot x\dot w_J H_w) = \orb_\CZ(g\dot x\orb_{\CZ_{\dot x}}(\dot w_J H_w)) = 
\orb_\CZ(g\dot xM^{w_J}).$$ 
The other assertions follow from Proposition~\ref{reduction}.
\end{proof}


\subsection{Main properties of $G^w$}

\begin{theorem}
\label{thm:zstablepieces1}
The $G^{w}$ for all $w\in \leftexp{I}{W}$ form a disjoint decomposition of~$G$.
\end{theorem}

\begin{proof}
We show this by induction on $\dim G$. In the base case $M=G$ we have $\leftexp{I}{W}=\{1\}$ and $H_1=G=G^1$ by Example~\ref{PequalGRem3}; hence the theorem is trivially true. Otherwise take an element $x\in\doubleexpIJ{W}$. By the induction hypothesis applied to the zip datum $E_{\CZ_{\dot x}}$ the subsets $M^{w_J}$ for $w_J\in \leftexp{I_x}{W_J}$ form a disjoint decomposition of~$M$. Thus by Proposition \ref{reduction} and Lemma \ref{lem:nicereduction2}, the subsets $G^{xw_J}$ for $w_J\in \leftexp{I_x}{W_J}$ form a disjoint decomposition of $Pg\dot xQ$. Combining this with the Bruhat decomposition \eqref{BruhatDecomp} it follows that the subsets $G^{xw_J}$ for all $x$ and $w_J$ form a disjoint decomposition of~$G$. But by Proposition \ref{howlettcor2} these are precisely the subsets $G^w$ for $w\in \leftexp{I}{W}$, as desired.
\end{proof}

\begin{theorem}
\label{GwVarDim}
For any $w \in \leftexp{I}{W}$ the subset $G^w$ is a nonsingular subvariety of $G$ of dimension $\dim P+\ell(w)$.
\end{theorem}

\begin{proof}
Again we proceed by induction on $\dim G$. If $M=G$, there is only one piece $G^1=G=P$ associated to $w=1$, and the assertion is clear. Otherwise write $w=xw_J$ as in Proposition~\ref{howlettcor2}. By the induction hypothesis the subset $M^{w_J}$ is a nonsingular subvariety of $M$ of dimension $\dim P_x+\ell(w_J)$. Thus by Propositions \ref{reduction} and \ref{reductiondimension1} and Lemma \ref{lem:nicereduction2} the subset $G^w$ is a nonsingular subvariety of dimension 
$$\bigl[\dim P_x + \ell(w_J)\bigr] + \dim P - \dim P_x + \ell(x)
\ =\ \dim P + \ell(x) + \ell(w_J).$$
By Proposition \ref{prop:howlett} the last expression is equal to $\dim P + \ell(w)$, as desired.
\end{proof}

\begin{theorem}
\label{thm:zstablepieces2}
For any $w \in \leftexp{I}{W}$, there is a bijection between the $H_w$-orbits $Y\subset H_w$ and the $E_\CZ$-orbits $X\subset G^w$, defined by $X=\orb_\CZ(g\dot wY)$ and satisfying
$$\codim(X\,{\subset}\,G^w) = \codim(Y\,{\subset}\,H_w).$$
\end{theorem}

\begin{proof}
If $M=G$, we have $w=1$ and $G=G^1=H_1$, and $E_\CZ\cong G$ acts on itself by the twisted conjugation $(h,h') \mapsto h\cdot h'\cdot\phi(h)^{-1}$. Thus the $E_\CZ$-orbits $X\subset G$ are precisely the cosets $gY$ for $H_1$-orbits $Y$ according to Definition \ref{def:HwMult}, which finishes that case.

If $M\not=G$ write $w=xw_J$ as in Proposition~\ref{howlettcor2}. Then $\CZ_{\dot x}$ has the frame $(M\cap B,T,1)$ by Proposition \ref{ReductionFrame}, and so by Lemma \ref{lem:iwreduction2} and the induction hypothesis there is a bijection between the $H_w$-orbits $Y\subset H_w$ and the $E_{\CZ_{\dot x}}$-orbits $X'\subset M^{w_J}$, defined by $X'=\orb_{\CZ_{\dot x}}(\dot w_JY)$ and satisfying $\codim(X'\,{\subset}\,M^{w_J}) = \codim(Y\,{\subset}\,H_w)$. By Proposition \ref{reduction} and Lemma \ref{lem:nicereduction2} there is a bijection between these $X'$ and the $E_\CZ$-orbits $X\subset G^w$, defined by $X=\orb_\CZ(g\dot xX')$. Moreover, since pushout and flat pullback preserve codimensions, the last statement in Proposition \ref{reduction} implies that 
$$\codim(X\,{\subset}\,G^w) 
= \codim\bigl(\orb_\CZ(g\dot xX') \,{\subset}\, \orb_\CZ(g\dot xM^{w_J})\bigr)
= \codim(X'\,{\subset}\,M^{w_J}).$$
Finally, since $\dot w = \dot x\dot w_J$ by \eqref{DotMultEq}, Lemma \ref{ozozw} shows that $X = \orb_\CZ(g\dot x\orb_{\CZ_{\dot x}}(\dot w_JY)) = \orb_\CZ(g\dot x\dot w_JY) \allowbreak = \orb_\CZ(g\dot wY)$, finishing the induction step.
\end{proof}


\subsection{Other descriptions of $G^w$}

\begin{lemma}\label{orbBB}
For any element $g'\in G$ we have
$$\orb_\CZ(gBg'B)\ =\ \orb_\CZ(gg'B)\ =\ \orb_\CZ(gBg').$$
\end{lemma}

\begin{proof}
Take any element $b\in B$. Then the condition \ref{FrameDef}~(b) implies that $p:=gbg^{-1}\in P$, and so there exists $q\in Q$ such that $(p,q)\in E_\CZ$. By the condition \ref{FrameDef}~(c) we then have $q\in B$. It follows that $gbg'B = pgg'Bq^{-1} \subset \orb_\CZ(gg'B)$. Since $b$ was arbitrary, this shows that $gBg'B \subset \orb_\CZ(gg'B)$, whence the first equality. A similar argument proves the second equality.
\end{proof}

\begin{theorem}
\label{Gw4} 
For any $w \in \leftexp{I}{W}$ we have
$$G^{w} 
= \orb_\CZ(g\dot wH_w) 
= \orb_\CZ(g\dot w(H_w\cap B))
= \orb_\CZ(g\dot wB)
= \orb_\CZ(gB\dot w)
= \orb_\CZ(gB\dot wB).$$
\end{theorem}

\begin{proof}
The first equation is Definition \ref{GwDef1} of~$G^w$, and the last two equations are cases of Lemma~\ref{orbBB}. The remaining two equations are proved by induction on $\dim G$. In the base case $M=G$ we have $w=1$ and $H_1=G$; hence the second term is $\orb_\CZ(gG) = G$, and the third and fourth terms are both equal to $\orb_\CZ(gB)$. By Proposition~\ref{lem:Brepr} applied to the isogeny $\phi\circ\inn(g)$ the latter is equal to~$G$, as desired.

In the case $M\not=G$ write $w=xw_J$ as in Proposition~\ref{howlettcor2}. Then $\CZ_{\dot x}$ has the frame $(M\cap B,T,1)$ by Proposition \ref{ReductionFrame}, and so by Lemma \ref{lem:iwreduction2} and the induction hypothesis we have
$$\orb_{\CZ_{\dot x}}(\dot w_JH_w) 
= \orb_{\CZ_{\dot x}}(\dot w_J(H_w\cap B))
= \orb_{\CZ_{\dot x}}(\dot w_J(M  \cap B)).$$
Using Lemma \ref{ozozw} this implies that
$$\orb_\CZ(g\dot x\dot w_JH_w) 
= \orb_\CZ(g\dot x\dot w_J(H_w\cap B))
= \orb_\CZ(g\dot x\dot w_J(M  \cap B)).$$

By \eqref{DotMultEq} we may replace $\dot x\dot w_J$ by $\dot w$ in these equations. Moreover, \eqref{ActionDef} and \eqref{ZipDecomp} show that $g\dot wB = g\dot w(M\cap B)V \subset \orb_\CZ(g\dot w(M  \cap B))$ and so $\orb_\CZ(g\dot wB) = \orb_\CZ(g\dot w(M\cap B))$. Thus both equations follow.
\end{proof}

\begin{example}
If $P$ is a Borel subgroup, then so is~$Q$, and we have $\leftexp{I}{W} = W$. The last equation in Theorem \ref{Gw4} then implies that $G^{w}=gB\dot wB$ for all $w \in W$.
\end{example}

For a further equivalent description of $G^w$ see Subsection~\ref{DiffCon}.


\section{Closure relation}
\label{Closure}

In this section, we determine the closure of $G^w$ in $G$ for any $w \in \leftexp{I}{W}$. To formulate a precise result recall that $\leq$ denotes the Bruhat order on~$W$. 

\begin{definition}
\label{RelDef}
For $w$, $w'\in\leftexp{I}{W}$ we write $w' \preccurlyeq w$ if and only if there exists $y \in W_{I}$ such that $yw'\psi(y)^{-1}\leq w$. 
\end{definition}

\begin{theorem}\label{thm:closure}
For any $w \in \leftexp{I}{W}$ we have
\begin{equation*}
\overline{G^{w}} = \coprod_{\substack{w^{\prime}\in \leftexp{I}{W}\\ w^{\prime}\preccurlyeq w}}G^{w^{\prime}}.
\end{equation*}
\end{theorem}

A direct consequence of this is:

\begin{corollary}
The relation $\preccurlyeq$ is a partial order on $\leftexp{I}{W}$.
\end{corollary}

\begin{remark}
The relation $\preccurlyeq$ has been introduced by He in \cite{He:GStablePieces} for a somewhat more special class of isomorphisms $\psi\colon W_I \iso W_J$. He gives a direct combinatorial proof that $\preccurlyeq$ is a partial order (Proposition~3.13 of loc.~cit.), which can be adapted to our more general setting (see \cite{specfzips}, Section 4).
\end{remark}


The rest of this section is devoted to proving Theorem~\ref{thm:closure}. We will exploit the fact that the closure relation for the Bruhat decomposition of $G$ is known. Namely, for any $w\in W$ we have by \cite{springer}, Proposition 8.5.5: 
\UseTheoremCounterForNextEquation
\begin{equation}\label{BruhatClos}
\overline{B\dot wB} = \coprod_{\substack{w'\in W\\ w'\leq w}}B\dot w'B.
\end{equation}

\begin{lemma}
\label{lem:closure}
For any $w \in W$ we have
\begin{equation*}
\overline{\orb_\CZ(gB\dot wB)} = \bigcup_{\substack{w' \in W \\ w'\leq w}}\orb_\CZ(g\dot w'B).
\end{equation*}
\end{lemma}

\begin{proof}
Let $B_\CZ\subset E_\CZ$ denote the subgroup of all elements $(u\ell,v\phi(\ell))$ with $u\in U$, $v\in V$, and $\ell\in L\cap \leftexp{g}{B}$. Then $E_\CZ/B_\CZ \cong L/(L\cap \leftexp{g}{B})$ is proper, and $gB\dot wB\subset G$ is a $B_\CZ$-invariant subvariety. Thus Lemma \ref{lem:properclosure} and \eqref{BruhatClos} imply that
$$\overline{\orb_\CZ(gB\dot wB)} \ =\ \orb_\CZ\bigl(\,\overline{gB\dot wB}\,\bigr)
\ = \bigcup_{w'\leq w}\orb_\CZ(gB\dot w'B).$$
The desired equality then follows from 
Lemma~\ref{orbBB}.
\end{proof}

\begin{lemma}
\label{lem:lem1}
For any $w$, $v \in W$ and $b \in B$ there exists $u \in W$ such that $u \leq v$ and $\dot w b \dot v \in B\dot{w}\dot{u}B$.
\end{lemma}

\begin{proof}
We prove the statement by induction on $\ell(v)$. If $v=1$, we may take $u=1$. For the induction step write $v=v^{\prime}s$ for some simple reflection $s$ such that $\ell(v^{\prime})=\ell(v)-1$. By the induction hypothesis there exists $u^{\prime}\leq v^{\prime}$ such that $\dot wb\dot v^{\prime}\in B\dot w \dot u^{\prime}B$. Hence $ \dot w b\dot v \in  B \dot w \dot u^{\prime} B \dot s\subset B\dot w \dot u^{\prime}\dot sB \cup B\dot w \dot u^{\prime}B$, so either $u=u^{\prime}s$ or $u=u^{\prime}$ will have the required property. 
\end{proof}

\begin{lemma}
\label{lem:kern}
For any $z\in W$ and $w\in \leftexp{I}{W}$ and $v\in W_I$ such that $z\leq w\psi(v)$, there exists $y\in W_I$ such that $yz\psi(y)^{-1}\leq vw$.
\end{lemma}

\begin{proof}
Choose reduced expressions for $w$ and $v$ as products of simple reflections. Since $\psi(I)=J$, this also yields a reduced expression for $\psi(v)$. Together this yields an expression for $w\psi(v)$ as a product of simple reflections, which is not necessarily reduced. However, by \cite{ccg}, Theorem 2.2.2 a reduced expression for $w\psi(v)$ can be obtained from the given one by possibly deleting some factors. By the definition of the Bruhat order, the assumption $z\leq w\psi(v)$ means that a reduced expression for $z$ is obtained from this by deleting further factors, if any. Let $y'$ denote the product of all factors remaining from~$w$. Since all factors in the reduced expression for $v$ lie in~$I$, the product of all factors remaining from~$\psi(v)$ is equal to $\psi(y)$ for some $y\in W_I$. By construction we then have $z=y'\psi(y)$, and so $yz\psi(y)^{-1} = yy'$. But the assumptions on $w$ and $v$ imply that $\ell(vw) = \ell(v)+\ell(w)$; hence the product of the given reduced expressions for $v$ and $w$ is a reduced expression for $vw$. By construction $yy'$ is obtained from that product by possibly deleting some factors, so we deduce that $yy'\leq vw$, as desired.
\end{proof}

\begin{lemma}
\label{lem:bw}
For any $w \in \leftexp{I}{W}$ and $w' \in W$ and $b$, $b^{\prime} \in B$ such that $\orb_\CZ(g\dot wb)=\orb_\CZ(g\dot w'b^{\prime})$ there exists $y \in W_{I}$ such that $yw\psi(y)^{-1}\leq w'$. 
\end{lemma}

\begin{proof}
We proceed by induction on $\dim G$. In the base case $M=G$ we have $w=1$ and may take $y=1$. So assume that $M\not=G$. Write $w=xw_J$ as in Proposition \ref{howlettcor2} with $x\in\doubleexp{I}{W}{J}$ and $w_J\in\leftexp{I_x}{W_J}$. From $\orb_\CZ(g\dot wb)=\orb_\CZ(g\dot w'b^{\prime})$ we deduce that $Pg\dot xQ = Pg\dot wQ = Pg\dot w'Q$, which in view of \eqref{BruhatDecomp} implies that $w' \in W_IxW_J$. Write $w'=v'xw_J'$ for $v'\in W_I$ and $w_J'\in\leftexp{I_x}{W_J}$, as in Proposition~\ref{prop:howlett}.

Recall that $\phi(g\dot v'g^{-1})\in\Norm_M(T)$ is a representative of $\psi(v')\in W_J$. Thus by Lemma \ref{lem:lem1}, there exists $u \in W$ such that $u \leq \psi(v')$ and $\dot x\dot w_J' b'\phi(g\dot v'g^{-1}) \in B\dot x\dot w_J'\dot{u}B$. The first condition implies that $u\in W_J$, the Weyl group of~$M$. The action of $E_\CZ$ and the second condition imply
$$\orb_\CZ(g\dot w'b')
\,=\, \orb_\CZ(g\dot v'\dot x\dot w_J'b')
\,=\, \orb_\CZ\bigl(g\dot x\dot w_J'b'\phi(g\dot v'g^{-1})\bigr)
\,\subset\, \orb_\CZ(gB\dot x\dot w_J'\dot{u}B).$$
Here the last term is equal to $\orb_\CZ(g\dot x\dot w_J'\dot{u}B)$ by Lemma~\ref{orbBB}. Thus there exists $b''\in B$ such that 
$$\orb_\CZ(g\dot x\dot w_Jb)
\,=\, \orb_\CZ(g\dot wb)
\,=\, \orb_\CZ(g\dot w'b')
\,=\, \orb_\CZ(g\dot x\dot w_J'\dot ub'').$$
By the action of $E_\CZ$ we may and do assume that $b$, $b''\in M\cap B$. Then $\dot w_Jb$ and $\dot w_J'\dot ub''$ lie in~$M$, and so Proposition \ref{reduction} implies that $\orb_{\CZ_{\dot x}}(\dot w_Jb) = \orb_{\CZ_{\dot x}}(\dot w_J'\dot ub'')$. By the induction hypothesis there therefore exists $y_x\in W_{I_x}$ such that 
$$y_xw_J\psi_x(y_x)^{-1} \leq w_J'u.$$
Now we work our way back up. Since both sides of the last relation lie in~$W_{J}$, and since $x\in W^{J}$, we deduce that
$$z\defeq xy_xw_J\psi_x(y_x)^{-1} \leq xw_J'u.$$
Recall that $u \leq \psi(v')$, which implies that $u=\psi(u')$ for some $u'\in W_I$ satisfying $u'\leq v'$. Also, note that $xw_J'\in \leftexp{I}{W}$ by Proposition~\ref{howlettcor2}. Thus by Lemma \ref{lem:kern} there exists $y'\in W_I$ such that 
$$y'z\psi(y')^{-1}\leq u'xw_J'.$$
As $u'$ and $v'$ lie in $W_I$, and $xw_J'\in \leftexp{I}{W}$, we deduce that 
$$y'z\psi(y')^{-1}\leq u'xw_J' \leq v' xw_J' = w'.$$
Finally, since $\psi_x=\psi\circ\inn(x)$, we have 
$$y'z\psi(y')^{-1} 
= y' xy_xw_J\psi(xy_xx^{-1})^{-1} \psi(y')^{-1}
= (y'xy_xx^{-1})xw_J\psi(y'xy_xx^{-1})^{-1}
= yw\psi(y)^{-1}$$
with $y := y'xy_xx^{-1} \in W_I$. Thus $yw\psi(y)^{-1}\leq w'$, as desired.
\end{proof}

\begin{lemma}\label{lem:ztdense}
For any $w \in \leftexp{I}{W}$, the set $\orb_\CZ(g\dot wT)$ is dense in $G^{w}$. 
\end{lemma}

\begin{proof}
Theorem \ref{thm:zstablepieces2} implies that $\orb_\CZ(g\dot wT) = \orb_\CZ(g\dot wY)$, where $Y\subset H_w$ is the orbit of $T$ under twisted conjugation by~$H_w$. But Proposition~\ref{lem:Brepr}~(b) asserts that $Y$ is dense in~$H_w$. Thus $\orb_\CZ(g\dot wY)$ is dense in $\orb_\CZ(\,\overline{g\dot wY}\,) = \orb_\CZ(g\dot wH_w) = G^w$, as desired.
\end{proof}

\begin{proof}[Proof of Theorem \ref{thm:closure}]
Consider $w^{\prime} \in \leftexp{I}{W}$ such that $G^{w^{\prime}}\cap\overline{G^{w}}\not=\emptyset$. Then by Theorem \ref{Gw4} and Lemma \ref{lem:closure} there exist $b$, $b^{\prime}\in B$ and $w'' \in W$ such that $w'' \leq w$ and $\orb_\CZ(gw^{\prime}b)=\orb_\CZ(gw''b^{\prime})$. Lemma \ref{lem:bw} then implies that $yw'\psi(y)^{-1}\leq w''$ for some $y \in W_{I}$. Together it follows that $yw'\psi(y)^{-1}\leq w$, and hence $w^{\prime} \preccurlyeq w$, proving ``$\subset$''.

To prove ``$\supset$'' consider $w^{\prime}\in \leftexp{I}{W}$ with $w^{\prime}\preccurlyeq w$. By definition there exists $y \in W_I$ such that $w'' := yw^{\prime}\psi(y)^{-1}\leq w$. Lemma \ref{lem:closure} and Theorem \ref{Gw4} then show that $\orb_\CZ(g\dot w'' T) \subset\overline{G^{w}}$. Therefore
$$\orb_\CZ(g\dot w^{\prime}T)
\, =\, \orb_\CZ\bigl(g\dot y \dot w^{\prime}T \phi(g\dot yg^{-1})^{-1}\bigr)
\, =\, \orb_\CZ\bigl(g\dot y \dot w^{\prime}\phi(g\dot yg^{-1})^{-1} T\bigr)
\, =\, \orb_\CZ(g\dot w'' T)
\, \subset\, \overline{G^{w}}.$$
With Lemma \ref{lem:ztdense} for $\orb_\CZ(g\dot w^{\prime}T)$ we conclude that $G^{w^{\prime}}\subset \overline{G^{w}}$, as desired.
\end{proof}


\section{Orbitally finite zip data}
\label{Frobenius}

\begin{proposition}\label{OrbFinEqu}
The following assertions are equivalent:
\begin{itemize}
\item[(a)] For any $w \in \leftexp{I}{W}$, the number of fixed points of the endomorphism $\phi_{\dot w}=\phi\circ\inn(g\dot w)$ of $H_w$ from Definition \ref{def:HwMult} is finite.
\item[(b)] For any $w \in \leftexp{I}{W}$ the $E_\CZ$-invariant subvariety $G^{w}$ is a single orbit under~$E_\CZ$.
\item[(c)] The number of orbits of $E_\CZ$ on $G$ is finite.
\end{itemize}
\end{proposition}

\begin{proof}
If (a) holds, the Lang-Steinberg Theorem \ref{thm:lang-steinberg} shows that the orbit of $1\in H_w$ under twisted conjugation is all of~$H_w$, and by Theorem \ref{thm:zstablepieces2} this implies~(b). The implication (b)$\implies$(c) is trivial. So assume~(c). Then again by Theorem \ref{thm:zstablepieces2}, the number of orbits in $H_w$ under twisted conjugation by $\phi_{\dot w}$ is finite for any $w \in \leftexp{I}{W}$. In particular there exists an open orbit; let $h$ be an element thereof. Then for dimension reasons its stabilizer is finite. But
$$\Stab_{H_w}(h)
\ =\ \bigl\{h'\in H_w\bigm|h'h\phi_{\dot w}(h')^{-1}=h\bigr\}
\ =\ \bigl\{h'\in H_w\bigm|h'=h\phi_{\dot w}(h')h^{-1}\bigr\}$$
is also the set of fixed points of the endomorphism $\inn(h)\circ\phi_{\dot w}$ of~$H_w$. Thus the Lang-Steinberg Theorem \ref{thm:lang-steinberg} implies that $\{h'h\phi_{\dot w}(h')^{-1}h^{-1}\mid h'\in H_w\} = H_w$. After right multiplication by $h$ this shows that the orbit of $h$ is all of~$H_w$. We may thus repeat the argument with the identity element in place of~$h$, and deduce that the set of fixed points of $\phi_{\dot w}$ on $H_w$ is finite, proving~(a).
\end{proof}

\begin{definition}
\label{def:frobenius}
We call $\CZ$ \emph{orbitally finite} if the conditions in Proposition \ref{OrbFinEqu} are met.
\end{definition}

\begin{proposition}
\label{dphi=0}
If the differential of $\phi$ at $1$ vanishes, then $\CZ$ is orbitally finite.
\end{proposition}

\begin{proof}
If the differential of $\phi$ vanishes, then so does the differential of $\phi_{\dot w}=\phi\circ\inn(g\dot w)|H_w$ for any $w \in \leftexp{I}{W}$. Let $H_w^f$ denote the fixed point locus of~$\phi_{\dot w}$, which is a closed algebraic subgroup. Then the restriction $\phi_{\dot w}|H_w^f$ is the identity and its differential is zero. This is possible only when $\dim H_w^f=0$, that is, when $H_w^f$ is finite.
\end{proof}

\begin{remark}
In particular Proposition \ref{dphi=0} applies when the base field has characteristic $p>0$ and the isogeny $\phi$ is a relative Frobenius $L\to L^{(p^r)}\cong M$.
\end{remark}

Since $g\dot w\in G^w$ by Definition \ref{GwDef1}, we can now rephrase condition \ref{OrbFinEqu} (b) and Theorems \ref{thm:zstablepieces1}, \ref{GwVarDim}, and \ref{thm:closure} as follows:

\begin{theorem}\label{thm:representatives}
Assume that $\CZ$ is orbitally finite. Then:
\begin{itemize}
\item[(a)] For any $w \in \leftexp{I}{W}$ we have $G^w=\orb_\CZ(g\dot w)$.
\item[(b)] The elements $g\dot w$ for $w \in \leftexp{I}{W}$ form a set of representatives for the $E_\CZ$-orbits in~$G$.
\item[(c)] For any $w \in \leftexp{I}{W}$ the orbit $\orb_\CZ(g\dot w)$ has dimension $\dim P+\ell(w)$.
\item[(d)] For any $w \in \leftexp{I}{W}$ the closure of $\orb_\CZ(g\dot w)$ is the union of $\orb_\CZ(g\dot w')$ for all $w' \in \leftexp{I}{W}$ with $w'\preccurlyeq w$.
\end{itemize}
\end{theorem}


\section{Point stabilizers}\label{Stabilizer}

In this section we study the stabilizer in $E_\CZ$ of an arbitrary element $g'\in G$. Take $w\in\leftexp{I}{W}$ such that $g' \in G^w$. Then Theorem \ref{thm:zstablepieces2} shows that $g'$ is conjugate to $g\dot wh$ for some $h\in H_w$. Thus it suffices to consider the stabilizer of $g\dot wh$. 

Recall from Definition \ref{def:HwMult} that $H_w$ acts on itself by twisted conjugation with the isogeny~$\phi_{\dot w}$, which is defined as the restriction of $\phi\circ\inn(g\dot w)$.

\begin{theorem}
\label{thm:stabilizer}
For any $w\in\leftexp{I}{W}$ and $h\in H_w$ the stabilizer $\Stab_{E_\CZ}(g\dot wh)$ is the semi-direct product of a connected unipotent normal subgroup with the subgroup
\UseTheoremCounterForNextEquation
\begin{equation}\label{OhDoch}
\bigl\{ \bigl(\inn(g\dot w)(h'),\phi(\inn(g\dot w)(h'))\bigr) \bigm| h'\in \Stab_{H_w}(h) \bigr\}.
\end{equation}
\end{theorem}

\begin{proof}
For any $h'\in H_w$ we have $\inn(g\dot w)(h') \cdot g\dot wh \cdot \phi(\inn(g\dot w)(h'))^{-1} = g\dot wh$ if and only if $h'h\phi(\inn(g\dot w)(h'))^{-1} = h$ if and only if $h'\in \Stab_{H_w}(h)$. This implies that \eqref{OhDoch} is a subgroup of $\Stab_{E_\CZ}(g\dot wh)$. 

For the rest we proceed by induction on $\dim G$. If $M=G$, we have $w=1$ and $\dot w=1$ and $G=H_1$, and $(g',\phi(g'))\in E_\CZ$ acts on $G$ by the twisted conjugation $g'' \mapsto g'g''\phi(g')^{-1}$. Under left translation by $g$ this corresponds to the action of $H_w$ on itself, so that $\Stab_{E_\CZ}(g\dot wh)$ is precisely the subgroup \eqref{OhDoch} and the normal subgroup is trivial.

If $M\not=G$ write $w=xw_J$ as in Proposition~\ref{howlettcor2}. Then $\dot w=\dot x\dot w_J$ and $\CZ_{\dot x}$ has the frame $(M\cap B,T,1)$, and Proposition \ref{stabilizerreduction} shows that $\Stab_{E_\CZ}(g\dot wh)$ is an extension of $\Stab_{E_{\CZ_{\dot x}}}(\dot w_Jh)$ by a connected unipotent normal subgroup. Moreover, by the induction hypothesis $\Stab_{E_{\CZ_{\dot x}}}(\dot w_Jh)$ is the semi-direct product of a connected unipotent normal subgroup with the subgroup
\UseTheoremCounterForNextEquation
\begin{equation}\label{OhNein}
\bigl\{ \bigl(\inn(\dot w_J)(h'),\phi_{\dot x}(\inn(\dot w_J)(h'))\bigr) \bigm| h'\in \Stab_{H_w}(h) \bigr\}.
\end{equation}
Furthermore a direct calculation shows that the projection in Proposition \ref{stabilizerreduction} sends the subgroup \eqref{OhDoch} isomorphically to the subgroup \eqref{OhNein}. Since any extension of connected unipotent groups is again connected unipotent, the theorem follows.
\end{proof}

\begin{remark}\label{EvensLu}
For the stabilizer, a similar result was obtained by Evens and Lu (\cite{EvensLu:LagrangianII} Theorem~3.13).
\end{remark}

If the differential of $\phi$ at $1$ vanishes, we can also describe the infinitesimal stabilizer in the Lie algebra. Since in that case the zip datum is orbitally finite by Proposition \ref{dphi=0}, it suffices to consider the stabilizer of $g\dot w$.

\begin{theorem}
\label{thm:liealgebradimension1}
Assume that the differential of $\phi$ at $1$ vanishes. For any $w\in \leftexp{I}{W}$ let $w=xw_J$ be the decomposition from Proposition \ref{howlettcor2}. Then the infinitesimal stabilizer of $g\dot w$ in the Lie algebra of $E_\CZ$ has dimension $\dim V-\ell(x)$. 
\end{theorem}

\begin{proof}
Since $d\phi=0$, we have $\Lie E_\CZ = \Lie P \times \Lie V \subset \Lie (P\times Q)$. Thus an arbitrary tangent vector of $E_\CZ$ at $1$ has the form $(1{+}dp,1{+}dv)$ for $dp\in\Lie P$ and $dv\in\Lie V$, viewed as infinitesimal elements of $P$ and $V$ in Leibniz's sense. That element stabilizes $g\dot w$ if and only if $(1{+}dp)g\dot w(1{+}dv)^{-1}=g\dot w$. This condition is equivalent to $dp \cdot g\dot w - g\dot w\cdot dv=0$, or again to $dp = \Ad_{g\dot w}(dv)$. The dimension is therefore $\dim(\Lie P \cap \Ad_{g\dot w}(\Lie V)) = \dim (\Lie \leftexp{g^{-1}}{P} \cap \Lie \leftexp{\dot w}{V})$. As both $\leftexp{g^{-1}}{P}$ and $\leftexp{\dot w}{V}$ are normalized by~$T$, the dimension is just the number of root spaces in the last intersection. This number is 
\begin{eqnarray*}
    \#\bigl[(\Phi^+\cup\Phi_I)\cap w(\Phi^+\setminus\Phi_J)\bigr]
&=& \#\bigl\{\alpha \in \Phi^+ \setminus \Phi_J \bigm| w\alpha \in \Phi^+ \cup \Phi_I\bigr\} \\
&=& \dim V - \#\bigl\{\alpha \in \Phi^+ \setminus \Phi_J \bigm| w\alpha \in \Phi^- \setminus \Phi_I\bigr\}.
\end{eqnarray*}
By Lemma \ref{lem:refinedlength} it is therefore $\dim V - \ell(x)$, as desired.
\end{proof}

\begin{remark}
\label{thm:liealgebradimension2}
The dimension in Theorem \ref{thm:liealgebradimension1} depends only on the first factor of $w=xw_J$ and thus only on the Bruhat cell $P g\dot wQ$. Since that Bruhat cell is an irreducible variety and in general composed of more than one $E_\CZ$-orbit, these orbits have different dimensions. Thus the corresponding point stabilizers in $E_\CZ$ have different dimension, while the dimension of their Lie algebra stabilizer is constant. Therefore the scheme-theoretic stabilizer of $g\dot w$ is in general not reduced.
\end{remark}


\section{Abstract zip data}\label{Abstract}

By Theorem \ref{thm:zstablepieces1} the subsets $G^{w}$ for all $w\in \leftexp{I}{W}$ form a disjoint decomposition of~$G$ satisfying ${g\dot w\in G^w}$. It is natural to ask which other elements of the form $g\dot w'$ for $w' \in W$ are contained in a given~$G^w$. When $\CZ$ is orbitally finite, by Theorem \ref{thm:representatives} this question is equivalent to asking which elements $g\dot w'$ for $w' \in W$ lie in the same $E_\CZ$-orbit. This problem turns out to depend only on the groups $W_I\subset W$ and the homomorphism~$\psi$ and can therefore be studied abstractly. We return to this situation at the end of this section.

\subsection{Abstract groups}

\begin{definition}
An \emph{abstract zip datum} is a tuple $\CA = (\Gamma,\Delta,\psi)$ consisting of a group $\Gamma$, a subgroup $\Delta$, and a homomorphism $\psi\colon\Delta\to\Gamma$.
\end{definition}

Fix such an abstract zip datum~$\CA$. For any $\gamma\in\Gamma$, the collection of subgroups $E$ of $\leftexp{\gamma^{-1}}{\!\Delta}$ satisfying $\psi\circ\inn(\gamma)(E) = E$ possesses a unique largest element, namely the subgroup generated by all such subgroups.

\begin{definition}
\label{Ew}
For any $\gamma\in\Gamma$ we let $E_\gamma$ denote the unique largest subgroup of $\leftexp{\gamma^{-1}}{\!\Delta}$ satisfying $\psi\circ\inn(\gamma)(E_\gamma) = E_\gamma$.
\end{definition}

\begin{lemma}
\label{graus}
For any $\gamma\in\Gamma$ and $\delta\in\Delta$ and $\epsilon\in E_\gamma$, we have $E_{\delta\gamma\epsilon\psi(\delta)^{-1}} = \leftexp{\psi(\delta)}{E_\gamma}$.
\end{lemma}

\begin{proof}
Abbreviate $\gamma' := \delta\gamma\epsilon\psi(\delta)^{-1}$. Then the calculation 
$\psi\bigl(\,\leftexp{\gamma'}{\!(}\leftexp{\psi(\delta)}{E_\gamma})\bigr)
=\psi(\,\leftexp{\delta\gamma\epsilon}{E_\gamma})
=\leftexp{\psi(\delta)}{\psi(\,\leftexp{\gamma}{E_\gamma})}
=\leftexp{\psi(\delta)}{E_\gamma}$
and the definition of $E_{\gamma'}$ imply that $\leftexp{\psi(\delta)}{E_\gamma} \subset E_{\gamma'}$. In particular, $\epsilon' := \leftexp{\psi(\delta)}{\epsilon}{}^{-1}$ is an element of~$E_{\gamma'}$. Since $\gamma = \delta'\gamma'\epsilon'\psi(\delta')^{-1}$ with $\delta' := \delta^{-1}\in\Delta$, a calculation like the first shows that $\leftexp{\psi(\delta')}{E_{\gamma'}} \subset E_\gamma$. Together it follows that $E_{\gamma'} = \leftexp{\psi(\delta)}{E_\gamma}$, as desired.
\end{proof}

\begin{definition}
\label{gnagna}
For any $\gamma$, $\gamma'\in\Gamma$ we write $\gamma'\sim\gamma$ if and only if there exist $\delta\in\Delta$ and $\epsilon\in E_\gamma$ such that $\gamma' = \delta\gamma\epsilon\psi(\delta)^{-1}$.
For any $\gamma\in\Gamma$ we abbreviate $\orb_\CA(\gamma) \defeq \{\gamma'\in\Gamma\mid\gamma'\sim\gamma\}$.
\end{definition}

\begin{lemma}
This is an equivalence relation.
\end{lemma}

\begin{proof}
Reflexivity is clear, and symmetry was shown already in the proof of Lemma~\ref{graus}. To prove transitivity, suppose that $\gamma' = \delta\gamma\epsilon\psi(\delta)^{-1}$ for $\delta\in\Delta$ and $\epsilon\in E_\gamma$ and $\gamma'' = \delta'\gamma'\epsilon'\psi(\delta')^{-1}$ for $\delta'\in\Delta$ and $\epsilon'\in E_{\gamma'}$. Then $\leftexp{\psi(\delta)^{-1}}{\!\epsilon}{}' \in E_\gamma$ by Lemma~\ref{graus}, and so $\gamma'' = \delta'\delta\gamma\epsilon\psi(\delta)^{-1}\epsilon'\psi(\delta')^{-1} = \delta''\gamma\epsilon''\psi(\delta'')^{-1}$ for $\delta'' := \delta'\delta\in\Delta$ and $\epsilon'' := \epsilon\;\leftexp{\psi(\delta)^{-1}}{\!\epsilon}{}' \in E_\gamma$, as desired.
\end{proof}

\begin{theorem}\label{eqclasscard} 
If $\Delta$ is finite, each equivalence class in $\Gamma$ has cardinality $\#\Delta$ and the number of equivalence classes is $[\Gamma:\Delta]$.
\end{theorem}

\begin{proof}
Take any $\gamma\in\Gamma$; then the group $E_\gamma \subset \leftexp{\gamma^{-1}}{\!\Delta}$ is finite, too.
Consider the surjective map ${\Delta\times E_\gamma\onto\orb_\CA(\gamma)}$, $(\delta,\epsilon)\mapsto \delta\gamma\epsilon\psi(\delta)^{-1}$. Two elements $(\delta,\epsilon)$, $(\delta',\epsilon') \in \Delta\times E_\gamma$ lie in the same fiber if and only if 
$\delta\gamma\epsilon\psi(\delta)^{-1} = \delta'\gamma\epsilon'\psi(\delta')^{-1}$
if and only if $\epsilon\psi(\delta^{-1}\delta') = \gamma^{-1}(\delta^{-1}\delta')\gamma\epsilon'$.
With $\epsilon'' := \gamma^{-1}(\delta^{-1}\delta')\gamma \allowbreak \in\leftexp{\gamma^{-1}}{\!\Delta}$ this is equivalent to $\epsilon\psi(\,\leftexp{\gamma}{\epsilon}{}'') = \epsilon''\epsilon'$.
Since $\epsilon$, $\epsilon'\in E_\gamma$, this equation implies that the subgroup generated by $E_\gamma$ and $\epsilon''$ is mapped onto itself under $\psi\circ\inn(\gamma)$. By maximality it is therefore equal to $E_\gamma$, and so $\epsilon''\in E_\gamma$. Together we find that the elements in the same fiber as $(\delta,\epsilon)$ are precisely the elements $(\delta',\epsilon')$ with $\delta' = \delta\,\leftexp{\gamma}{\epsilon}{}''$ and $\epsilon' = (\epsilon'')^{-1}\epsilon\psi(\,\leftexp{\gamma}{\epsilon}{}'')$ for some $\epsilon''\in E_\gamma$. Thus each fiber has cardinality $\#E_\gamma$, and so the image has cardinality $\#\Delta$, proving the first assertion. The second assertion is a direct consequence of the first.
\end{proof}

We can also perform an induction step as in Section~\ref{Bruhat} for abstract zip data, obtaining analogues of Lemma \ref{lem:iwreduction2} and Proposition~\ref{reduction}. For this fix an element $\xi\in\Gamma$, say in a set of representatives for the double quotient $\Delta\backslash\Gamma/\psi(\Delta)$. Then Definitions \ref{Ew} and \ref{gnagna} imply that the equivalence class of any $\gamma\in\Delta\xi\psi(\Delta)$ is again contained in $\Delta\xi\psi(\Delta)$. 

\begin{construction}\label{abscons}
Set $\Gamma_\xi \defeq \psi(\Delta)$ and $\Delta_\xi \defeq \psi(\Delta)\cap\leftexp{\xi^{-1}}{\!\Delta}$, and let $\psi_\xi\colon \Delta_\xi\to\Gamma_\xi$ denote the restriction of $\psi\circ\inn(\xi)$. This defines a new, possibly smaller, abstract zip datum
$$\CA_\xi\defeq (\Gamma_\xi,\Delta_\xi,\psi_\xi).$$
\end{construction}

\begin{lemma}\label{absiwreduction}
For any $\gamma\in\Gamma_\xi$, the group $E_{\xi\gamma}$ associated by Definition \ref{Ew} to the pair $(\CA,\xi\gamma)$ is equal to the group associated to the pair $(\CA_\xi,\gamma)$.
\end{lemma}

\begin{proof}
Since $\gamma\in\Gamma_\xi=\psi(\Delta)$, Definition~\ref{Ew} implies that
$$E_{\xi\gamma}\ \subset\ \psi(\Delta)\cap\leftexp{\gamma^{-1}\xi^{-1}}{\!\Delta}
\ =\ \leftexp{\gamma^{-1}}{\bigl(}\psi(\Delta)\cap\leftexp{\xi^{-1}}{\!\Delta}\bigr)
\ =\ \leftexp{\gamma^{-1}}{\!\Delta_\xi}$$
and that $E_{\xi\gamma} = \psi\circ\inn(\xi\gamma)(E_{\xi\gamma}) = \psi_\xi\circ\inn(\gamma)(E_{\xi\gamma})$. Since $E_{\xi\gamma}$ is the largest subgroup of $\leftexp{\gamma^{-1}\xi^{-1}}{\!\Delta}$ with this property, it is also the largest in $\leftexp{\gamma^{-1}}{\!\Delta_\xi}$.
\end{proof}

\begin{proposition}\label{absabsreduction}
There is a bijection between $\CA_\xi$-equivalence classes in $\Gamma_\xi$ and $\CA$-equiva\-lence classes in $\Delta\xi\psi(\Delta)$, defined by $\orb_{\CA_\xi}(\gamma)\mapsto \orb_\CA(\xi\gamma)$ and $\orb_{\CA_\xi}(\gamma) = \Gamma_\xi\cap \xi^{-1}\!\orb_\CA(\xi\gamma)$.
\end{proposition}

\begin{proof}
Take any $\gamma$, $\gamma'\in \Gamma_\xi$. Then $\gamma'\in \xi^{-1}\!\orb_\CA(\xi\gamma)$ if and only if $\xi\gamma' = \delta\xi\gamma\epsilon\psi(\delta)^{-1}$ for some $\delta\in\Delta$ and $\epsilon\in E_{\xi\gamma}$. Writing $\delta=\leftexp{\xi}{\delta'}$ this is equivalent to $\gamma' = \delta'\gamma\epsilon\psi(\,\leftexp{\xi}{\delta'})^{-1}$ for $\delta'\in\leftexp{\xi^{-1}}{\!\Delta}$ and $\epsilon\in E_{\xi\gamma}$. In this equation $\gamma$ and $\gamma'$ and $\psi(\,\leftexp{\xi}{\delta'})$ lie in $\Gamma_\xi=\psi(\Delta)$ by assumption, and so does $\epsilon\in E_{\xi\gamma} \subset \psi(\Delta)$ by Definition~\ref{Ew}. Thus the equation requires that $\delta'$ lies in~$\psi(\Delta)$, and so a fortiori in $\psi(\Delta)\cap \leftexp{\xi^{-1}}{\!\Delta} = \Delta_\xi$. In view of Lemma \ref{absiwreduction} the condition is thus equivalent to $\gamma' \in \orb_{\CA_\xi}(\gamma)$, proving the equation at the end of the proposition.

That equation implies that the map $\orb_{\CA_\xi}(\gamma)\mapsto \orb_\CA(\xi\gamma)$ from $\CA_\xi$-equivalence classes in $\Gamma_\xi$ to $\CA$-equivalence classes in $\Delta\xi\psi(\Delta)$ is well-defined and injective. But any element of $\Delta\xi\psi(\Delta)$ has the form $\delta\xi\gamma$ for $\delta\in\Delta$ and $\gamma\in\Gamma_\xi$ and is therefore equivalent to $\xi\gamma\psi(\delta) \in \xi\Gamma_\xi$. Thus the map is also surjective, and we are done.
\end{proof}


\subsection{Coxeter groups}
\label{wuerg}

\begin{definition}
Let $W$ be a Coxeter group with a finite set of simple reflections~$S$. Let $\psi\colon W_I\iso W_J\subset W$ be an isomorphism of Coxeter groups with $\psi(I)=J$ for subsets $I$, $J\subset S$. Then $\CA:=(W,W_I,\psi)$ is an abstract zip datum that we call \emph{of Coxeter type}.
\end{definition}

Fix such an abstract zip datum of Coxeter type~$\CA$. Recall that $\doubleexp{I}{W}{J}$ is a set of representatives for the double quotient $W_I\backslash W/W_J$. We will apply the induction step from Proposition \ref{absabsreduction} to $x\in\doubleexp{I}{W}{J}$. As in Proposition \ref{nicereduction1} set $I_x\defeq J\cap\leftexp{x^{-1}}{\!I}$ and $J_x\defeq \psi(I\cap\leftexp{x}{J})$, which are both subsets of~$J$. Then $W_J=\psi(W_I)$, and $W_{I_x} = \psi(W_I) \cap \leftexp{x^{-1}}{W_I}$ by Proposition~\ref{thm:kilmoyer}, and $\psi_x :=\psi\circ\inn(x)$ induces an isomorphism $\psi_x\colon W_{I_x}\iso W_{J_x}$ such that $\psi_x(I_x)=J_x$. Thus the new abstract zip datum from Construction \ref{abscons} is $\CA_x := (W_J, W_{I_x},\psi_x)$ and hence again of Coxeter type. Using this we obtain the following analogue of Theorem~\ref{thm:zstablepieces1}, which also has been previously proved by He (\cite{He:MinimalLengthCoxeter} Corollary~2.6).

\begin{theorem}
\label{Abszstablepieces1}
For $\CA$ of Coxeter type $\leftexp{I}{W}$ is a set of representatives for the equivalence classes in~$W$.
\end{theorem}

\begin{proof}
We prove this by induction on~$\#S$. If $I=S$, we have $W_I=W_J=W$ and so $E_w=W$ for every $w\in W$. Then there is exactly one equivalence class, represented by the unique element of $\leftexp{I}{W}=\{1\}$, and the assertion holds.

Otherwise we have $\#I<\#S$. Take any $x\in\doubleexp{I}{W}{J}$. Then by the induction hypothesis $\leftexp{I_x}{W_J}$ is a set of representatives for the $\CA_x$-equivalence classes in~$W_J$. Thus Proposition \ref{absabsreduction} implies that $x\,\leftexp{I_x}{W_J}$ is a set of representatives for the $\CA$-equivalence classes in $W_IxW_J$. Varying $x$, Proposition \ref{howlettcor2} implies that $\leftexp{I}{W}$ is a set of representatives for the equivalence classes in~$W$, as desired.
\end{proof}

For use in Section \ref{Diff} we include the following results.

\begin{lemma}\label{Bij1}
\begin{itemize}
\item[(a)] For any $w\in\leftexp{I}{W}$ there exists $y\in W_I$ such that $w'\defeq yw\psi(y)^{-1} \in W^J$.
\item[(b)] The element $w'$ in (a) is independent of~$y$.
\end{itemize}
\end{lemma}

\begin{proof}
For (a) we use induction on~$\#S$. If $I=S$, we have $\leftexp{I}{W}=\{1\}$ and $w=1$, and so $y=1$ does the job. Otherwise $\#I<\#S$. Write $w=xw_J$ as in Proposition \ref{howlettcor2} with $x\in\doubleexp{I}{W}{J}$ and $w_J\in\leftexp{I_x}{W_J}$. Then by the induction hypothesis applied to $\CA_x$ there exists $y'\in W_{I_x}$ such that 
$$w_J' \ \defeq\ y'w_J\psi_x(y')^{-1} 
       \ \in\ W_J^{J_x}
       \ =\ W_{\psi(I)}^{\psi(I\cap\,\leftexp{x}{\!J})}
       \ =\ \psi\bigl(W_I^{I\cap\,\leftexp{x}{\!J}}\,\bigr).$$
Setting $y \defeq\psi^{-1}(y'w_J)\in W_I$ and using the definition of $\psi_x$ we deduce that
\begin{eqnarray*}
w'\ \defeq\ yw\psi(y)^{-1} &=& \psi^{-1}(y'w_J)\cdot xw_J\cdot(y'w_J)^{-1} \\
               &=& \psi^{-1}(y'w_J)\cdot xy^{\prime-1}x^{-1}\cdot x \\
               &=& \psi^{-1}\bigl( y'w_J\psi_x(y')^{-1} \bigr) \cdot x \\
               &=& \psi^{-1}(w_J') \cdot x \ \in\ W_I^{I\cap\,\leftexp{x}{\!J}}\cdot x.
\end{eqnarray*}
By Proposition \ref{howlettcor2dual} the right hand side is contained in~$W^J$, showing~(a). 

To prove (b) consider another element $y'\in W_I$ such that $w''\defeq y'w\psi(y')^{-1} \in W^J$. Then with $\tilde y := \psi(y'y^{-1}) \in W_J$ we have $w''= y'y^{-1}w'\psi(y)\psi(y')^{-1} = \psi^{-1}(\tilde y) w' \tilde y^{-1}$ and hence $w^{\prime\prime-1} = \tilde y w^{\prime-1} \psi^{-1}(\tilde y)^{-1}$. Now observe that on replacing $(I,J,\psi)$ by $(J,I,\psi^{-1})$ we obtain another abstract zip datum $\CA' \defeq(W,W_J,\psi^{-1})$ dual to~$\CA$. The last equality then shows that $w^{\prime\prime-1}$ and $w^{\prime-1}$ are equivalent according to Definition \ref{gnagna} for~$\CA'$. Since these elements also lie in $\leftexp{J}{W}$, Theorem \ref{Abszstablepieces1} applied to $\CA'$ shows that they are equal. Therefore $w''=w'$, as desired. 
\end{proof}

\begin{proposition}\label{Bij}
There exists a unique bijection $\sigma\colon\leftexp{I}{W}\to W^J$ with the property that for any $w\in\leftexp{I}{W}$ there exists $y\in W_I$ such that $\sigma(w) = yw\psi(y)^{-1}$.
\end{proposition}

\begin{proof}
The existence of a unique \emph{map} $\sigma\colon\leftexp{I}{W}\to W^J$ with the stated property is equivalent to Lemma~\ref{Bij1}. By applying the same lemma to the abstract zip datum $\CA' \defeq(W,W_J,\psi^{-1})$ in place of $\CA$ we find that for any $w'\in\leftexp{J}{W}$ there exists $y'\in W_J$ such that $w\defeq y'w'\psi^{-1}(y')^{-1} \in W^I$, and the element $w$ is independent of~$y'$. After replacing $(w',w)$ by $(w^{\prime-1},w^{-1})$ this means that for any $w'\in W^J$ there exists $y'\in W_J$ such that $w\defeq \psi^{-1}(y')w'y^{\prime-1} \in \leftexp{I}{W}$, and the element $w$ is independent of~$y'$. But with $y := \psi^{-1}(y')^{-1} \in W_I$ the last equation is equivalent to $w' = yw\psi(y)^{-1}$, and so for any $w'\in W^J$ there exists a unique $w\in\leftexp{I}{W}$ with $w'=\sigma(w)$. In other words the map is bijective, as desired.
\end{proof}

\begin{proposition}\label{BijLength}
The bijection in Proposition~\ref{Bij} satisfies $\ell(w)=\ell(\sigma(w))$ for all $w\in\leftexp{I}{W}$.
\end{proposition}

\begin{proof}
Write the defining relation in the form $yw=\sigma(w)\psi(y)$. Here $y\in W_I$ and $w\in\leftexp{I}{W}$ imply that $\ell(yw) = \ell(y)+\ell(w)$, and similarly $\sigma(w)\in W^J$ and $\psi(y)\in W_J$ imply that $\ell(\sigma(w)\psi(y)) = \ell(\sigma(w))+\ell(\psi(y))$. Moreover, since $\psi$ sends simple reflections to simple reflections, it satisfies $\ell(\psi(y))=\ell(y)$. Together it follows that $\ell(w)=\ell(\sigma(w))$.
\end{proof}

\begin{lemma}\label{Bij2}
Let $\sigma\colon\leftexp{I}{W}\to W^J$ be the bijection from Proposition~\ref{Bij}. For any $x\in\doubleexp{I}{W}{J}$ let $\sigma_x\colon\leftexp{I_x}{W_J}\to W^{J_x}_J$ denote the bijection obtained by applying Proposition~\ref{Bij} to~$\CA_x$. Then for all $w_J\in\leftexp{I_x}{W_J}$ we have $\sigma(xw_J) = \psi^{-1}(\sigma_x(w_J))\cdot x$.
\end{lemma}

\begin{proof}
The proof of Lemma \ref{Bij1}~(a) shows that $\sigma(w) = w' = \psi^{-1}(w_J') \cdot x$ where $w_J'=\sigma_x(w_J)$, as desired.
\end{proof}

\begin{remark}\label{CompareHe}
Propositions~\ref{Bij} and~\ref{BijLength} can also be deduced from more general results of He (\cite{He:MinimalLengthCoxeter}~Proposition~4.3).
\end{remark}


\subsection{Back to algebraic groups}

Now we return to the situation and the notations of the preceding sections. Clearly the connected algebraic zip datum~$\CZ$ gives rise to an abstract zip datum of Coxeter type $\CA \defeq (W,W_I,\psi)$, which by Proposition \ref{Indep1} is independent of the frame, up to unique isomorphism. Theorem \ref{thm:zstablepieces1} implies that for any $w'\in W$ the element $g\dot w'$ lies in $G^w$ for a unique ${w\in\leftexp{I}{W}}$.

\begin{theorem}
\label{gut}
For any $w'\in W$ and $w\in\leftexp{I}{W}$ we have $g\dot w' \in G^w$ if and only if $w'\sim w$ with respect to~$\CA$.
\end{theorem}

\begin{proof}
We prove this by induction on~$\#S$. If $J=S$, there is exactly one $G^w$ for $w=1$ and exactly one $\CA$-equivalence class in~$W$, so the assertion holds. Otherwise we have $\#J<\#S$. Write $w=xw_J$ with $x\in\doubleexp{I}{W}{J}$ and $w_J\in\leftexp{I_x}{W_J}$, as in Proposition~\ref{howlettcor2}. Then by \eqref{BruhatDecomp} and Lemma \ref{lem:nicereduction2} the condition $g\dot w'\in G^w$ requires that $w'\in W_IxW_J$, and so does the condition $w'\sim w$ by the remarks in Subsection~\ref{wuerg}. It therefore suffices to consider $w'=yxw_J'$ with $y\in W_I$ and $w_J'\in\leftexp{I_x}{W_J}$, as in Proposition~\ref{prop:howlett}. But then $w'\sim xw_J'\psi(y)$ with respect to~$\CA$, and $g\dot w' = g\dot y\dot x\dot w_J'$ is in the same $E_\CZ$-orbit as $g\dot x\dot w_J'\phi(\leftexp{g}{\dot y})$. After replacing $w'$ by $xw_J'\psi(y)$ we may thus assume that $w'=xw_J'$ for some $w_J'\in W_J$. Then Proposition \ref{reduction} and Lemma \ref{lem:nicereduction2} show that $g\dot w'\in G^w$ if and only if $\dot w_J'\in M^{w_J}$. By the induction hypothesis this is equivalent to $w_J'\sim w_J$ with respect to~$\CA_x$. By Proposition \ref{absabsreduction} this in turn is equivalent to $w'\sim w$ with respect to~$\CA$, as desired.
\end{proof}

Combining Theorems \ref{thm:representatives} and \ref{gut} we deduce:

\begin{corollary}
If $\CZ$ is orbitally finite, then for any $w$, $w'\in W$ the elements $g\dot w$ and $g\dot w'$ lie in the same $E_\CZ$-orbit if and only if $w\sim w'$ with respect to~$\CA$.
\end{corollary}


\section{Non-connected algebraic zip data}\label{nonconnected}

In this section we generalize the main results of Sections \ref{Gw} and \ref{Closure} to non-connected groups. Throughout we denote a not necessarily connected linear algebraic group by~$\hat G$, its identity component by~$G$, and its finite group of connected components by $\pi_0(\hat G) := \hat G/G$; and similarly for other letters of the alphabet. Note that the unipotent radical $R_uG$ is a normal subgroup of~$\hat G$. Any homomorphism $\hat\phi\colon\hat G\to\hat H$ restricts to a homomorphism $\phi\colon G\to H$.

\begin{definition}\label{NCZipDatumDef}
An \emph{algebraic zip datum} is a tuple $\hat\CZ = (\hat G,\hat P,\hat Q,\hat\varphi)$ consisting of a linear algebraic group $\hat G$ with subgroups $\hat P$ and $\hat Q$ and a homomorphism $\hat\varphi\colon\hat P/\CR_uP\to\hat Q/\CR_uQ$, such that $\CZ := (G,P,Q,\phi)$ is a connected algebraic zip datum. The zip group $E_{\hat\CZ} \subset \hat P\times\hat Q$, its action on~$\hat G$, and the orbit $\orb_{\hat\CZ}(X)$ of a subset $X\subset\hat G$ are defined in exact analogy to \eqref{ZipGroupDef}, \eqref{ActionDef}, and~\eqref{OrbDef}. 
\end{definition}

Throughout this section we fix an algebraic zip datum $\hat\CZ = (\hat G,\hat P,\hat Q,\hat\varphi)$ with associated connected algebraic zip datum $\CZ = (G,P,Q,\phi)$. We fix a frame $(B,T,g)$ of $\CZ$ and use the other pertaining notations from Sections \ref{AlgZipData} through~\ref{Gw}. We also define 
$$\hat W\defeq \Norm_{\hat G}(T)/T\qquad\hbox{and}\qquad\Omega\defeq (\Norm_{\hat G}(B) \cap \Norm_{\hat G}(T))/T,$$
so that $\Omega\cong \pi_0(\hat G)$ and  $\hat W= W\rtimes\Omega$. For each $\omega\in\Omega$ we fix a representative $\dot\omega\in\Norm_{\Ghat}(B)\cap\Norm_{\Ghat}(T)$, and for $\hat w=w\omega\in\hat W$ with $w\in W$ and $\omega\in\Omega$ we set $\dot{\hat w} \defeq \dot w\dot\omega\in\Norm_{\Ghat}(T)$.

Note that by definition $E_\CZ$ is the identity component of~$E_{\hat\CZ}$. Thus to study the $E_{\hat\CZ}$-orbits in~$\hat G$, we first study the orbits under $E_\CZ$ and then the action of $E_{\hat\CZ}/E_\CZ$ on them.


\begin{lemma}
\label{ReductionToConnected}
For any $\omega\in\Omega$ the \emph{conjugate} connected algebraic zip datum
$$\leftexp{\dot\omega}{\!\CZ} \defeq (G,P,\leftexp{\dot\omega}{Q},\inn(\dot\omega)\circ\phi)$$
has zip group $E_{\,\leftexp{\dot\omega}{\!\CZ}} = \{(p,\leftexp{\dot\omega}{q})\mid(p,q)\in E_\CZ\}$ and frame $(B,T,g)$, and the isomorphism of varieties $G\to G\dot\omega$, $g'\mapsto g'\dot\omega$ induces a bijection from the $E_{\,\leftexp{\dot\omega}{\!\CZ}}$-orbits in $G$ to the $E_\CZ$-orbits in~$G\dot\omega$.
\end{lemma}

\begin{proof}
Direct calculation.
\end{proof}

\begin{lemma}
\label{ConnectedOrbits}
The subsets $\orb_\CZ(gB\dot{\hat w}B)$ for all $\hat w\in \leftexp{I}W\Omega$ form a disjoint decomposition of~$\hat G$.
\end{lemma}

\begin{proof}
Take any $\omega\in\Omega$. Then by Theorems \ref{thm:zstablepieces1} and \ref{Gw4} the subsets $\orb_{\,\leftexp{\dot\omega}{\!\CZ}}(gB\dot wB)$ for all $w\in \leftexp{I}{W}$ form a disjoint decomposition of~$G$. Thus by Lemma \ref{ReductionToConnected} the subsets $\orb_{\CZ}(gB\dot wB\dot\omega)$ for all $w\in \leftexp{I}{W}$ form a disjoint decomposition of~$G\dot\omega$. Since $\dot\omega\in\Norm_{\hat G}(B)$ by assumption, the latter subset is equal to $\orb_{\CZ}(gB\dot w\dot\omega B)$. By varying $\omega$ the proposition follows.
\end{proof}


Next define $\hat L\defeq \Norm_{\hat P}(L)$ and $\hat M\defeq \Norm_{\hat Q}(M)$, so that $\hat P=U\rtimes\hat L$ and $\hat Q=V\rtimes\hat M$, and $\hat\phi$ can be identified with a homomorphism $\hat L\to\hat M$. Set
$$\begin{tabular}{ll}
\llap{$\hat W_I$}$\;\defeq \Norm_{\leftexp{g^{-1}}{\!\hat L}}(T)/T$,\qquad\qquad &
\llap{$\Omega_I$}$\;\defeq (\Norm_{\leftexp{g^{-1}}{\!\hat L}}(B) \cap \Norm_{\leftexp{g^{-1}}{\!\hat L}}(T))/T$, \\
\llap{$\hat W_J$}$\;\defeq \Norm_{\hat M}(T)/T$,\qquad\qquad &
\llap{$\Omega_J$}$\;\defeq (\Norm_{\hat M}(B) \cap \Norm_{\hat M}(T))/T$. \\
\end{tabular}$$
These groups are subgroups of $\hat W$ and satisfy
$$\begin{tabular}{ll}
\llap{$\hat W_I$}$\;= W_I\rtimes\Omega_I$,\qquad\qquad &
\llap{$\Omega_I$}$\;\cong \pi_0(\hat L)\cong \pi_0(\hat P)$, \\
\llap{$\hat W_J$}$\;= W_J\rtimes\Omega_J$,\qquad\qquad &
\llap{$\Omega_J$}$\;\cong \pi_0(\hat M)\cong \pi_0(\hat Q)$. \\
\end{tabular}$$
Also $\hat\phi\circ\inn(g)$ induces a homomorphism $\hat\psi\colon \hat W_I\to\hat W_J$ extending $\psi\colon W_I\to W_J$ and sending $\Omega_I$ to~$\Omega_J$.
Moreover, the elements $(\leftexp{g}{\dot\omega},\hat\phi(\leftexp{g}{\dot\omega}))$ for all $\omega\in\Omega_I$ are representatives of the connected components of~$E_{\hat\CZ}$.

\begin{lemma}
\label{ConnectedOrbitsAction}
\begin{itemize}
\item[(a)] The map $(\upsilon,\hat w)\mapsto\upsilon\hat w\hat\psi(\upsilon)^{-1}$ defines a left action of $\Omega_I$ on $\leftexp{I}W\Omega$.
\item[(b)] Take any $\upsilon\in\Omega_I$ and $\hat w\in \leftexp{I}W\Omega$ and abbreviate $\hat w' \defeq \upsilon\hat w\hat\psi(\upsilon)^{-1} \in \leftexp{I}W\Omega$. Then the element $(\leftexp{g}{\dot\upsilon},\hat\phi(\leftexp{g}{\dot\upsilon})) \in E_{\hat\CZ}$ sends $\orb_\CZ(gB\dot{\hat w}B)$ to $\orb_\CZ(gB\dot{\hat w}{}'B)$.
\end{itemize}
\end{lemma}

\begin{proof}
Conjugation by $\Omega_I$ preserves the set of simple reflections $I$ and thus the subset ${\leftexp{I}{W}\subset W}$. In (a) we therefore have $\upsilon\hat w\hat\psi(\upsilon)^{-1} = \leftexp{\upsilon}{\hat w}\cdot\upsilon\hat\psi(\upsilon)^{-1} \in \leftexp{I}{W}\Omega\cdot\Omega = \leftexp{I}{W}\Omega$, as desired. In (b) the elements $\dot\upsilon$ and $\hat\phi(\leftexp{g}{\dot\upsilon})$ normalize~$B$; hence the image is 
$$\leftexp{g}{\dot\upsilon} \cdot \orb_\CZ(gB\dot{\hat w}B) \cdot \hat\phi(\leftexp{g}{\dot\upsilon})^{-1}
\ =\ \orb_\CZ\bigl( \leftexp{g}{\dot\upsilon} gB\dot{\hat w}B\hat\phi(\leftexp{g}{\dot\upsilon})^{-1} \bigr)
\ =\ \orb_\CZ\bigl( gB\dot\upsilon\dot{\hat w}\hat\phi(\leftexp{g}{\dot\upsilon})^{-1}B \bigr).$$
As $\dot\upsilon\dot{\hat w}\hat\phi(\leftexp{g}{\dot\upsilon})^{-1}$ differs from $\dot{\hat w}{}'$ by an element of~$T$, this proves~(b).
\end{proof}


For any $\hat w \in \leftexp{I}W\Omega$ we now define
\UseTheoremCounterForNextEquation
\begin{equation}\label{DefineGwEqGeneral}
\hat G^{\hat w}\defeq\orb_{\hat\CZ}(gB\dot{\hat w} B),
\end{equation}
which is independent of the representative~$\dot{\hat w}$. Lemma \ref{ConnectedOrbitsAction} implies that $\hat G^{\hat w}$ is the union of $\orb_\CZ(gB\dot{\hat w}{}'B)$ for all $\hat w'$ in the $\Omega_I$-orbit of $\hat w$ under the action in \ref{ConnectedOrbitsAction}~(a). Thus $\hat G^{\hat w}$ depends only on $\hat w$ modulo~$\Omega_I$, and with Lemma \ref{ConnectedOrbits} we conclude:

\begin{theorem}
\label{NCGw}
The subsets $\hat G^{\hat w}$ for all $\hat w\in \leftexp{I}{W}\Omega$ modulo the action of $\Omega_I$ from \ref{ConnectedOrbitsAction}~(a) form a disjoint decomposition of~$\hat G$.
\end{theorem}


To describe the closure relation between the subsets $\hat G^{\hat w}$ we define analogues of the Bruhat order $\leq$ on $\hat W=W\Omega$ and of the relation $\preccurlyeq$ from Definition \ref{RelDef} on $\leftexp{I}W\Omega$:

\begin{definition}
\label{NCPO1}
For $\hat w=w\omega$ and $\hat w'=w'\omega'$ with $w$, $w'\in W$ and $\omega$, ${\omega'\in\Omega}$ we write $\hat w'\leq \hat w$ if and only if $w' \leq w$ and $\omega'=\omega$.
\end{definition}

\begin{definition}
\label{NCPO2}
For $\hat w$, $\hat w'\in\leftexp{I}{W}\Omega$ we write $\hat w'\preccurlyeq\hat w$ if and only if there exists $\hat y\in\hat W_I$ such that $\hat y\hat w'\hat\psi(\hat y)^{-1}\leq \hat w$. 
\end{definition}

\begin{theorem}\label{NCclosure}
For any $\hat w \in \leftexp{I}{W}\Omega$ we have
\begin{equation*}
\overline{\hat G^{\hat w}}\ = 
\bigcup_{\substack{\hat w'\in \leftexp{I}{W}\Omega\\ \hat w'\preccurlyeq\hat w}} \hat G^{\hat w'}.
\end{equation*}
\end{theorem}

\begin{proof}
Write $\hat w=w\omega$ with $w\in \leftexp{I}{W}$ and $\omega\in\Omega$. Then the conjugate zip datum $\leftexp{\dot\omega}{\!\CZ}$ has the isogeny $\inn(\dot\omega)\circ\phi\colon L\to \leftexp{\dot\omega}{M}$ and hence the induced isomorphism of Weyl groups ${\inn(\omega)\circ\psi\colon} W_I \iso \leftexp{\omega}{W_J} = W_{\leftexp{\omega}{\!J}}$. Thus Theorems \ref{Gw4} and \ref{thm:closure} and Definition \ref{RelDef} imply that
$$\strut\hskip3cm\llap{$\overline{\orb_{\leftexp{\dot\omega}{\!\CZ}}(gB\dot wB)}$}
\ =\ \bigcup_{w'} \rlap{$\displaystyle\orb_{\leftexp{\dot\omega}{\!\CZ}}(gB\dot w'B),$}\hskip3cm\strut$$
where the union ranges over all $w'\in \leftexp{I}{W}$ such that $yw'\:\leftexp{\omega}{\psi}(y)^{-1}\leq w$ for some $y \in W_{I}$. Note that this inequality is equivalent to $yw'\omega\psi(y)^{-1}\leq w\omega$ by Definition \ref{NCPO1}. Thus with Lemma \ref{ReductionToConnected} we deduce that 
$$\strut\hskip5cm\llap{$\overline{\orb_{\CZ}(gB\dot{\hat w} B)}
\ =\ \overline{\orb_{\CZ}(gB\dot w\dot\omega B)}$}
\ =\ \bigcup_{w'}\rlap{$\displaystyle \orb_{\CZ}(gB\dot w'\dot\omega B)
\ =\ \bigcup_{\hat w'} \orb_{\CZ}(gB\dot{\hat w}' B),$}\hskip5cm\strut$$
where the last union ranges over all $\hat w'\in \leftexp{I}{W}\Omega$ such that $y\hat w'\psi(y)^{-1}\leq \hat w$ for some $y \in W_{I}$. By taking the union of conjugates of this under $(\leftexp{g}{\dot\upsilon},\hat\phi(\leftexp{g}{\dot\upsilon})) \in E_{\hat\CZ}$ for all $\upsilon\in\Omega_I$ we obtain the closure of~$\hat G^{\hat w}$. By Lemma \ref{ConnectedOrbitsAction} the right hand side then yields the union of $\orb_\CZ(gB\dot{\hat w}{}''B)$ for all $\hat w'' = \upsilon\hat w'\hat\psi(\upsilon)^{-1}$ with $y\hat w'\psi(y)^{-1}\leq \hat w$ for some $\upsilon\in\Omega_I$ and $y \in W_{I}$. But here $\hat y := y\upsilon^{-1}$ runs through the group $W_I\Omega_I = \hat W_I$ and the inequality is equivalent to
$$\hat y\hat w''\hat\psi(\hat y)^{-1} \ =\ y\upsilon^{-1} \hat w'\hat\psi(\upsilon)\psi(y)^{-1}\ \leq\hat w.$$
By Definition \ref{NCPO2} these $\hat w''$ are precisely the elements of $\leftexp{I}{W}\Omega$ satisfying $\hat w''\preccurlyeq\hat w$.
\end{proof}


Finally, let us call $\hat\CZ$ \emph{orbitally finite} if the conjugates $\leftexp{\dot\omega}{\!\CZ}$ are orbitally finite for all ${\omega\in\Omega}$. This holds in particular when the differential of $\hat\phi$ at $1$ vanishes, because then 
we can apply Proposition~\ref{dphi=0} to~$\leftexp{\dot\omega}{\!\CZ}$. Combining 
Theorem \ref{thm:representatives} with the remarks leading up to Theorem \ref{NCGw} we deduce:

\begin{theorem}\label{HatOrbFin}
Assume that $\hat\CZ$ is orbitally finite. Then:
\begin{itemize}
\item[(a)] For any $\hat w \in \leftexp{I}{W}\Omega$ we have $\hat G^{\hat w}=\orb_{\hat\CZ}(g\dot{\hat w})$.
\item[(b)] If $\hat w \in \leftexp{I}{W}\Omega$ runs through a system of representatives for the action of $\Omega_I$ from \ref{ConnectedOrbitsAction}~(a), then $g\dot{\hat w}$ runs through a set of representatives for the $E_{\hat\CZ}$-orbits in~$\hat G$.
\end{itemize}
\end{theorem}


\section{Dual parametrization}
\label{Diff}

The decomposition of $G$ from Theorem \ref{thm:zstablepieces1} is parametrized in a natural way by elements of~$\leftexp{I}{W}$. In this section we translate that parametrization into an equally natural parametrization by elements of~$W^J$, which was used by Lusztig and He (see Section~\ref{CompareLusztig}). We also carry out the corresponding translation in the non-connected case.


\subsection{The connected case}
\label{DiffCon}

For any $w \in W^J$ we set
\UseTheoremCounterForNextEquation
\begin{equation}\label{GwDef2}
G^{w} \defeq \orb_\CZ\bigl(gB\dot wB\bigr).
\end{equation}
Note that this does not depend on the representative $\dot w$ of $w$ and conforms to Definition \ref{GwDef1} by Theorem~\ref{Gw4}. In Proposition \ref{Bij} we have already established a natural bijection $\sigma\colon\leftexp{I}{W}\to W^J$.

\begin{theorem}\label{dualparam}
For any $w \in \leftexp{I}{W}$ we have $G^w=G^{\sigma(w)}$.
\end{theorem}

\begin{proof}
If $I=J=S$, we have $\leftexp{I}{W}=W^J=\{1\}$ and so $w=\sigma(w)=1$; hence the assertion holds trivially. Otherwise $\#I<\#S$. Write $w=xw_J$ as in Proposition \ref{howlettcor2} with $x\in\doubleexp{I}{W}{J}$ and $w_J\in\leftexp{I_x}{W_J}$, and let $\sigma_x\colon\leftexp{I_x}{W_J}\to W^{J_x}_J$ denote the bijection obtained by applying Proposition~\ref{Bij} to~$\CA_x$. Then $\CZ_{\dot x}$ has the frame $(M\cap B,T,1)$ by Proposition \ref{ReductionFrame}, and so the induction hypothesis implies that
$$M^{w_J}\ =\ M^{\sigma_x(w_J)}\ =\ \orb_{\CZ_{\dot x}}\bigl((M\cap B)\dot\sigma_x(w_J)(M\cap B)\bigr).$$
By Lemma \ref{orbBB} this is equal to $\orb_{\CZ_{\dot x}}((M\cap B)\dot\sigma_x(w_J))$, and so by Lemmas \ref{lem:nicereduction2} and \ref{ozozw} we have
$$G^w \ =\ \orb_\CZ(g\dot xM^{w_J})
 \ =\ \orb_\CZ\bigl(g\dot x\orb_{\CZ_{\dot x}}((M\cap B)\dot\sigma_x(w_J))\bigr)
 \ =\ \orb_\CZ\bigl(g\dot x(M\cap B)\dot\sigma_x(w_J)\bigr).$$
Recall from Lemma \ref{Bij2} that $\sigma(w) = w_Ix$ with $w_I\defeq\psi^{-1}(\sigma_x(w_J))\in W_I$. It follows that $\sigma_x(w_J) = \psi(w_I)$ and therefore $\dot\sigma_x(w_J) \in T\cdot\phi(\leftexp{g}{\dot w_I})$ and $\dot\sigma(w) \in T\cdot\dot w_I\dot x$. Since $T\subset M\cap B$, using the action \eqref{ActionDef} of $E_\CZ$ we deduce that 
$$G^w\ =\ \orb_\CZ\bigl(g\dot x(M\cap B)\phi(\leftexp{g}{\dot w_I})\bigr)
\ =\ \orb_\CZ\bigl(g\dot w_I\dot x(M\cap B)\bigr)
\ =\ \orb_\CZ\bigl(g\dot\sigma(w)(M\cap B)\bigr).$$
Using \eqref{ActionDef} and \eqref{ZipDecomp} for the action of $V$, respectively Lemma \ref{orbBB}, we conclude that 
$$G^w\ =\ \orb_\CZ\bigl(g\dot\sigma(w)B\bigr)
\ =\ \orb_\CZ\bigl(gB\dot\sigma(w)B\bigr)
\ =\ G^{\sigma(w)},$$
as desired.
\end{proof}

\begin{theorem}
\label{Dualzstablepieces1}
The $G^{w}$ for all $w\in W^J$ form a disjoint decomposition of~$G$ by nonsingular subvarieties of dimension $\dim P+\ell(w)$.
\end{theorem}

\begin{proof}
Combine Theorems \ref{thm:zstablepieces1}, \ref{GwVarDim}, \ref{dualparam} and Proposition~\ref{BijLength}.
\end{proof}


Next, in analogy to Definition \ref{RelDef} we define:

\begin{definition}
\label{RelDefJ}
For $w$, $w'\in W^J$ we write $w' \preccurlyeq w$ if and only if there exists $y \in W_{I}$ such that $yw'\psi(y)^{-1}\leq w$. 
\end{definition}

\begin{theorem}\label{dualclosure}
For any $w \in W^J$ we have
\begin{equation*}
\overline{G^{w}} = \coprod_{\substack{w^{\prime}\in W^J\\ w^{\prime}\preccurlyeq w}}G^{w^{\prime}}.
\end{equation*}
\end{theorem}

\begin{proof}
By combining Theorems \ref{dualparam} and \ref{thm:closure} we already know that $\overline{G^{w}}$ is the disjoint union of $G^{w'}$ for certain $w'\in W^J$; it only remains to determine which. 

First consider $w' \in W^J$ with $G^{w'}\subset\overline{G^{w}}$. Then $g\dot w'\in\overline{G^{w}}$, and so by Lemma \ref{lem:closure} there exist $b\in B$ and $w'' \in W$ such that $w'' \leq w$ and $\orb_\CZ(g\dot w')=\orb_\CZ(g\dot w''b)$. Set $\tilde w'\defeq\sigma^{-1}(w')\in\leftexp{I}{W}$ and take $y\in W_I$ satisfying $w'=y\tilde w'\psi(y)^{-1}$. Then $\dot w'=\dot y\dot{\tilde w}'t\phi(\leftexp{g}{\dot y})^{-1}$ for some $t\in T$, and thus $\orb_\CZ(gw')=\orb_\CZ(g\tilde w't)$. Therefore $\orb_\CZ(g\tilde w't)=\orb_\CZ(gw''b)$, and so Lemma \ref{lem:bw} implies that $y'\tilde w'\psi(y')^{-1}\leq w''$ for some $y'\in W_{I}$. Together it follows that 
$$(y'y^{-1})w'\psi(y'y^{-1})^{-1} = y'\tilde w'\psi(y')^{-1}\leq w''\leq w$$ 
and hence $w^{\prime} \preccurlyeq w$, proving ``$\subset$''.

Conversely consider $w' \in W^J$ with $w'\preccurlyeq w$, and take $y \in W_I$ such that $w'' := yw'\psi(y)^{-1}\allowbreak{\leq w}$. Lemma \ref{lem:closure} then shows that $\orb_\CZ(g\dot w'' T) \subset\overline{G^{w}}$. Therefore
$$\orb_\CZ(g\dot w^{\prime}T)
\, =\, \orb_\CZ\bigl(g\dot y \dot w^{\prime}T \phi(g\dot yg^{-1})^{-1}\bigr)
\, =\, \orb_\CZ\bigl(g\dot y \dot w^{\prime}\phi(g\dot yg^{-1})^{-1} T\bigr)
\, =\, \orb_\CZ(g\dot w'' T)
\, \subset\, \overline{G^{w}}.$$
Since also $\orb_\CZ(g\dot w^{\prime}T)\subset G^{w'}$, this with the preliminary remark on $\overline{G^{w}}$ shows that $G^{w'}\subset \overline{G^{w}}$, proving~``$\supset$''.
\end{proof}

\begin{remark}
\label{AnyW}
In Definitions \ref{GwDef1} and \ref{GwDef2} we have introduced the subsets $G^w\defeq \orb_\CZ\bigl(gB\dot wB\bigr)$ only for $w\in \leftexp{I}{W}\cup W^J$, not for arbitrary $w\in W$. Our results do not say anything directly about the latter. Note that in case $\phi$ is an isomorphism their closures have been determined in \cite{LuYakimov:GroupCompactification}~Theorem~5.2 and \cite{He:MinimalLengthCoxeter}~Proposition~5.8.
\end{remark}


\subsection{The non-connected case}
\label{DiffNonCon}

Now we return to the notations from Section~\ref{nonconnected}. We begin with an analogue of Proposition~\ref{Bij}: 

\begin{proposition}\label{HatBij}
There exists a unique bijection $\hat\sigma\colon\leftexp{I}{W}\Omega\to \Omega W^J$ with the property that for any $\hat w\in\leftexp{I}{W}\Omega$ there exists $y\in W_I$ such that $\hat\sigma(\hat w) = y\hat w\psi(y)^{-1}$.
\end{proposition}

\begin{proof}
The equation requires that $\hat\sigma(\hat w)\in\hat W$ lie in the same $W$-coset as~$\hat w$. Thus for any fixed $\omega\in\Omega$, we need a unique bijection $\leftexp{I}{W}\omega\to\omega W^J$ sending $w\omega$ to an element of the form $yw\omega\psi(y)^{-1}$ for some $y\in W_I$. Multiplying both elements on the right by $\omega^{-1}$ this amounts to a unique bijection $\leftexp{I}{W}\to \omega W^J\omega^{-1} = W^{\,\leftexp{\omega}{\!J}}$ sending $w$ to an element of the form $yw\omega\psi(y)^{-1}\omega^{-1}$
for some $y\in W_I$. But $\inn(\omega)\circ\psi\colon W_I\to \omega W_J\omega^{-1} = W_{\leftexp{\omega}{\!J}}$ is precisely the isomorphism associated to the conjugate connected algebraic zip datum $\leftexp{\dot\omega}{\!\CZ}$ from Lemma~\ref{ReductionToConnected}. Thus a unique bijection with that property exists by Proposition \ref{Bij} applied to $\leftexp{\dot\omega}{\!\CZ}$.
\end{proof}

For any $\hat w \in \Omega W^J$ we now define
\UseTheoremCounterForNextEquation
\begin{equation}\label{DefineGwEqGeneralDual}
\hat G^{\hat w}\defeq\orb_{\hat\CZ}(gB\dot{\hat w} B).
\end{equation}
Again this does not depend on the representative $\dot{\hat w}$ of $\hat w$ and conforms to Definition~\eqref{DefineGwEqGeneral}. 

\begin{theorem}\label{hatdualparam}
For any $\hat w \in \leftexp{I}{W}\Omega$ we have $\hat G^{\hat w}=\hat G^{\hat\sigma(\hat w)}$.
\end{theorem}

\begin{proof}
Write $\hat w = w\omega$ with $w\in\leftexp{I}{W}$ and $\omega\in\Omega$. In the proof of Proposition \ref{HatBij} we have seen that $\hat\sigma(\hat w) = w'\omega$, where $w'\in W^{\,\leftexp{\omega}{\!J}}$ is the image of $w$ under the isomorphism given by Proposition \ref{Bij} applied to $\leftexp{\dot\omega}{\!\CZ}$. Thus by Theorem \ref{dualparam} we  have $\orb_{\,\leftexp{\dot\omega}{\!\CZ}}(gB\dot wB) = \orb_{\,\leftexp{\dot\omega}{\!\CZ}}(gB\dot w'B)$ inside~$G$. On multiplying on the right by $\dot\omega$ and applying Lemma \ref{ReductionToConnected} to both sides we deduce that
$$\orb_{\CZ}(gB\dot{\hat w}B) = 
\orb_{\CZ}(gB\dot wB\dot\omega) = \orb_{\CZ}(gB\dot w'B\dot\omega)
= \orb_{\CZ}(gB\dot{\hat\sigma}(\hat w)B).$$
The desired equality follows from this by applying~$\orb_{\hat\CZ}$.
\end{proof}

\begin{lemma}
\label{DualConnectedOrbitsAction}
\begin{itemize}
\item[(a)] The map $(\upsilon,\hat w)\mapsto\upsilon\hat w\hat\psi(\upsilon)^{-1}$ defines a left action of $\Omega_I$ on $\Omega W^J$.
\item[(b)] The bijection $\hat\sigma\colon\leftexp{I}{W}\Omega\to \Omega W^J$ from Proposition \ref{HatBij} is $\Omega_I$-equivariant.
\end{itemize}
\end{lemma}

\begin{proof}
Take $\upsilon\in\Omega_I$ and $\hat w\in \Omega W^J$. To prove (a) observe that conjugation by $\hat\psi(\upsilon)\in\Omega_J$ preserves the set of simple reflections $J$ and thus the subset ${W^J\subset W}$. We therefore have $\upsilon\hat w\hat\psi(\upsilon)^{-1} = \upsilon\hat\psi(\upsilon)^{-1} \cdot \leftexp{\hat\psi(\upsilon)}{\hat w} \in \Omega\cdot\Omega W^J = \Omega W^J$, as desired. In (b) write $\hat\sigma(\hat w) = y\hat w\psi(y)^{-1}$ with $y\in W_I$. Then 
$$\upsilon\hat\sigma(\hat w)\hat\psi(\upsilon)^{-1} 
= (\upsilon y \upsilon^{-1}) (\upsilon\hat w\hat\psi(\upsilon)^{-1})\hat\psi(\upsilon y\upsilon^{-1})^{-1}
= \hat\sigma\bigl(\upsilon\hat w\hat\psi(\upsilon)^{-1}),$$
because the left hand side is in $\Omega W^I$ and $\upsilon y \upsilon^{-1}\in W_I$. This proves~(b).
\end{proof}

\begin{theorem}
\label{DualNCGw}
The subsets $\hat G^{\hat w}$ for all $\hat w\in\Omega W^J$ modulo the action of $\Omega_I$ from \ref{DualConnectedOrbitsAction}~(a) form a disjoint decomposition of~$\hat G$.
\end{theorem}

\begin{proof}
Combine Theorems \ref{NCGw} and \ref{hatdualparam} with Lemma \ref{DualConnectedOrbitsAction}.
\end{proof}

\begin{definition}
\label{NCPO2Dual}
For $\hat w$, $\hat w'\in\Omega W^J$ we write $\hat w'\preccurlyeq\hat w$ if and only if there exists $\hat y\in\hat W_I$ such that $\hat y\hat w'\hat\psi(\hat y)^{-1}\leq \hat w$. 
\end{definition}

\begin{theorem}\label{NCclosureDual}
For any $\hat w \in\Omega W^J$ we have
\begin{equation*}
\overline{\hat G^{\hat w}}\ = 
\bigcup_{\substack{\hat w'\in\Omega W^J\\ \hat w'\preccurlyeq\hat w}} \hat G^{\hat w'}.
\end{equation*}
\end{theorem}

\begin{proof}
Write $\hat w = w\omega$ with $\omega\in\Omega$ and $w\in W^{\,\leftexp{\omega}{\!J}}$. Applying Theorem \ref{dualclosure} to the conjugate zip datum $\leftexp{\dot\omega}{\!\CZ}$ shows that $\overline{\orb_{\,\leftexp{\dot\omega}{\!\CZ}}(gB\dot wB)}$ is the union of the subsets $\orb_{\,\leftexp{\dot\omega}{\!\CZ}}(gB\dot w'B)$ for all $w'\in W^{\,\leftexp{\omega}{\!J}}$ such that $yw'\omega\psi(y)^{-1}\omega^{-1}\leq w$ for some $y \in W_{I}$. On multiplying on the right by $\dot\omega$ and applying Lemma \ref{ReductionToConnected} to everything we deduce that
$\overline{\orb_{\CZ}(gB\dot{\hat w}B)} = \overline{\orb_{\CZ}(gB\dot wB\dot\omega)}$ is the union of the subsets $\orb_{\CZ}(gB\dot w'B\dot\omega) = \orb_{\CZ}(gB\dot w'\dot\omega B)$ for the same elements~$w'$. Writing $\hat w' = w'\omega$ this is equal to the union of the subsets $\orb_{\CZ}(gB\dot{\hat w}'B)$ for all $\hat w'\in\Omega W^J$ such that $y\hat w'\psi(y)^{-1}\leq w$ for some $y \in W_{I}$. The theorem follows from this by applying~$\orb_{\hat\CZ}$.
\end{proof}


\section{Generalization of certain varieties of Lusztig}\label{CompareLusztig}

In this section we consider a certain type of algebraic variety with an action of a reductive group $G$ whose orbit structure is closely related to the structure of the $E_\CZ$-orbits in $G$ for an algebraic zip datum~$\CZ$. Special cases of such varieties have been defined by Lusztig (\cite{Lusztig:parsheaves2}) and by Moonen and the second author in~\cite{moonwed}.


\subsection{The coset variety of an algebraic zip datum}
\label{cosetvar}

\begin{remark}\label{cosetvarNC}
To keep notations simple, we restrict ourselves to connected zip data, although everything in this section directly extends to non-connected ones by putting $\hat{\phantom{x}}$ in the appropriate places.
\end{remark}

In this section we use only the definition of algebraic zip data and the action of the associated zip group from Section~\ref{AlgZipData}, but none of the other theory or notations from the preceding sections, not even the concept of a frame. Fix a connected algebraic zip datum $\CZ = (G,P,Q,\varphi)$. Recall that $E_\CZ$ is a subgroup of $P\times Q$ and hence of $G\times G$. We also consider the image of $ G$ under the diagonal embedding $\Delta\colon G\hookrightarrow G\times G$, $g\mapsto(g,g)$. We are interested in the left quotient $\Delta(G)\backslash(G\times G)$ and the right quotient $(G\times G)/E_\CZ$.

The first is isomorphic to $G$ via the projection morphism 
\UseTheoremCounterForNextEquation
\begin{equation}\label{lambda}
\lambda\colon G\times G \to G,\ (g,h)\mapsto g^{-1} h.
\end{equation}
Turn the right action of $E_\CZ$ on $G\times G$ into a left action by letting $(p,q)\in E_\CZ$ act by right translation with $(p,q)^{-1}$. Then with $E_\CZ$ acting on $G$ as in the definition of algebraic zip data, a direct calculation 
shows that $\lambda$ is $E_\CZ$-equivariant.

To describe the second quotient recall that $\phi$ is a homomorphism $P/U\to Q/V$, where $U$ and $V$ denote the unipotent radicals of $P$ and~$Q$. Consider a left $P$-coset $X\subset G$ and a left $Q$-coset $Y\subset G$. Then $X/U$ is a right torsor over $P/U$, and $Y/V$ is a right torsor over $Q/V$. By a $P/U$-equivariant morphism $\Phi\colon X/U\to Y/V$ we mean a morphism satisfying $\Phi(\bar x\bar p) = \Phi(\bar x)\phi\bigl(\bar p)$ for all $\bar x\in X/U$ and $\bar p\in P/U$. 

\begin{definition}\label{cosetdef}
The \emph{coset space of $\CZ$} is the set $C_\CZ$ of all triples $(X,Y,\Phi)$ consisting of a left $P$-coset $X\subset G$, a left $Q$-coset $Y\subset G$, and a $P/U$-equivariant morphism $\Phi\colon X/U\to Y/V$.
\end{definition}

For any $X$, $Y$ as above and any $(g,h)\in G\times G$, left multiplication by $g$ induces an isomorphism $\ell_g\colon X/U\iso g X/U$, and left multiplication by $h$ induces an isomorphism $\ell_h\colon Y/V\iso h Y/V$. Therefore $(X,Y,\Phi) \mapsto \bigl(gX,hY,\ell_h\circ\Phi\circ\ell_g^{-1}\bigr)$ defines a left action of $G\times G$ on~$C_\CZ$. By applying this action to the canonical base point $(P,Q,\phi)\in C_\CZ$ we obtain a morphism
\UseTheoremCounterForNextEquation
\begin{equation}\label{rho}
\rho\colon G\times G\to C_\CZ,\ 
  (g,h) \mapsto 
  \bigl(gP,h Q,\ell_{h}\circ\phi\circ\ell_{g}^{-1}\bigr).
\end{equation}
Clearly this morphism is equivariant under the left action of $G\times G$ and hence under the subgroup $\Delta(G)$.

\begin{lemma}\label{quot}
There is a unique structure of algebraic variety on $C_\CZ$ such that
$\rho$ identifies $C_{\CZ}$ with the quotient variety $(G\times G)/E_\CZ$.
\end{lemma}

\begin{proof}
The action of $G\times G$ is obviously transitive on the set of all pairs $(X,Y)$. Moreover, any $P/U$-equivariant morphism of right torsors $P/U\to Q/V$ has the form $\bar p \mapsto \pi_{Q}(q) \phi(\bar p) = \ell_{q}\circ \phi(\bar p)$ for some $q\in Q$. Thus the subgroup $1\times Q$ acts transitively on the set of all triples of the form $(P,Q,\Phi)$. Together it follows that the action of $G\times G$ on $C_\CZ$ is transitive.

On the other hand $(g,h)$ lies in the stabilizer of $(P,Q,\phi)$ if and only if $g\in P$ and $h\in Q$ and $\ell_h\circ\phi\circ\ell_g^{-1} = \phi$. But under the first two of these conditions, we have for all $\bar p\in P/U$
$$\ell_{h}\circ\phi\circ\ell_{g}^{-1} (\bar p)
= \pi_{Q}(h) \phi\bigl( \pi_{P}(g)^{-1} \bar p\bigr)
= \pi_{Q}(h) \phi\bigl( \pi_{P}(g)\bigr)^{-1} \phi(\bar p),$$
and so the third condition is equivalent to $\phi\bigl( \pi_{P}(g)\bigr) = \pi_{Q}(h)$. Together this means precisely that $(g,h) \in E_\CZ$, which is therefore the stabilizer of $(P,Q,\phi)$.

It follows that $\rho$ induces a bijection $(G\times G)/E_\CZ \iso C_{\CZ}$. Since the quotient variety exists by \cite{efa}, Section 3.2, this yields the unique structure of algebraic variety on~$C_\CZ$.
\end{proof}

Following Lemma \ref{quot} we call $C_\CZ$ also the \emph{coset variety of~$\CZ$}. Recall from \cite{efa} Prop.\ 2.5.3 that the quotient of an algebraic group by an algebraic subgroup is always a torsor. To summarize we have therefore constructed morphisms with the following properties:
\UseTheoremCounterForNextEquation
\begin{equation}\label{summary}
\fbox{$\displaystyle
\begin{tabular}{r}
$E_\CZ$-equivariant\\[3pt]
$\Delta(G)$-torsor
\end{tabular}
\left\{\vcenter{\hsize=0pt\vrule width0pt height35pt}\right.
\vcenter{\xymatrix@C-15pt{
& G\times G \ar[dl]_(.6){\textstyle\mathstrut\lambda} \ar[dr]^(.6){\textstyle\mathstrut\rho} & \\
G && C_\CZ \\}}
\left.\vcenter{\hsize=0pt\vrule width0pt height35pt}\right\}
\begin{tabular}{l}
$\Delta(G)$-equivariant\\[3pt]
$E_\CZ$-torsor
\end{tabular}
$}
\end{equation}

Recall that the actions of $\Delta(G)$ and $E_\CZ$ on $G\times G$ commute and thus combine to an action of $\Delta(G)\times E_\CZ$. Therefore \eqref{summary} directly implies:

\begin{theorem}\label{stacks}
There are natural isomorphisms of algebraic stacks
$$\xymatrix{[E_\CZ\backslash G] & 
\bigl[(\Delta(G)\times E_\CZ)\backslash(G\times G)\bigr] 
\ar[l]_-{[\lambda]}^-{\sim} \ar[r]^-{[\rho]}_-{\sim} &
[\Delta(G)\backslash C_\CZ]. \\}$$
\end{theorem}

Even without stacks, we can deduce:

\begin{theorem}\label{orbitbijection}
\begin{itemize}
\item [(a)] There is a closure-preserving bijection between $E_\CZ$-invariant subsets $A\subset G$ and $\Delta(G)$-invariant subsets $B\subset C_\CZ$, defined by $A = \lambda(\rho^{-1}(B))$ and $B = \rho(\lambda^{-1}(A))$. 
\item [(b)] The subset $A$ in (a) is a subvariety, resp.\ a nonsingular subvariety, if and only if $B$ is one. In that case we also have $\dim A = \dim B$.
\item [(c)] In particular (a) induces a bijection between $E_\CZ$-orbits in $G$ and $\Delta(G)$-orbits in~$C_\CZ$.
\item [(d)] For any $g\in G$ and $(X,Y,\Phi)\in C_\CZ$ whose orbits correspond, there is an isomorphism
$$\Stab_{E_\CZ}(g)\ \cong\ \Stab_{\Delta(G)}((X,Y,\Phi)).$$
\end{itemize}
\end{theorem}

\begin{proof}
By \eqref{summary} any $\Delta(G)\times E_\CZ$-invariant subset of $G\times G$ must be simultaneously of the form $\lambda^{-1}(A)$ for an $E_\CZ$-invariant subset $A\subset G$ and of the form $\rho^{-1}(B)$ for a $\Delta(G)$-invariant subset $B\subset C_\CZ$. Then $A = \lambda(\rho^{-1}(B))$ and $B = \rho(\lambda^{-1}(A))$, giving the bijection in~(a). The bijection preserves closures because $\lambda$ and $\rho$ are smooth. This proves~(a), the first sentence in~(b), and the special case~(c). In (b) it also proves that $\dim A+\dim G = \dim B+\dim E_\CZ$. But $\dim G = 2\dim U+\dim L = 2\dim V+\dim M$ and $\dim L=\dim M$ imply that $\dim U=\dim V$, and thus using \eqref{ZipDecomp} that $\dim E_\CZ = \dim U + \dim L + \dim V =\dim G$. Therefore $\dim A = \dim B$, proving the rest of~(b).

In (c) by assumption there exists a point $\underline{x}\in G\times G$ such that $\lambda(\underline{x})$ lies in the $E_\CZ$-orbit of~$g$ and $\rho(\underline{x})$ lies in the $\Delta(G)$-orbit of $(X,Y,\Phi)$. Thus after replacing $\underline{x}$ by a suitable translate under $\Delta(G)\times E_\CZ$ we may assume that $\lambda(\underline{x})=g$ and $\rho(\underline{x})=(X,Y,\Phi)$. Then the fact that $\lambda$ and $\rho$ are torsors implies that the two projection morphisms
$$\qquad\xymatrix{{\mathstrut}\Stab_{E_\CZ}(g) & 
{\mathstrut}\Stab_{\Delta(G)\times E_\CZ}(\underline{x}) \ar[l] \ar[r] &
{\mathstrut}\Stab_{\Delta(G)}((X,Y,\Phi)) \\}$$
are isomorphisms, proving~(c).
(The isomorphism may depend on the choice of $\underline{x}$.)
\end{proof}

With Theorem \ref{orbitbijection} we can translate many results about the $E_\CZ$-action on $G$ from the preceding sections to the $\Delta(G)$-action on~$C_\CZ$, in particular Theorems \ref{thm:zstablepieces1}, \ref{GwVarDim}, \ref{thm:closure}, \ref{thm:representatives}, \ref{thm:stabilizer}, and their counterparts from Sections \ref{nonconnected} and~\ref{Diff}.


\subsection{Algebraic zip data associated to an isogeny of~$G$}
\label{AZDisog}

In this subsection we consider algebraic zip data whose isogeny extends to an isogeny on all of~$G$. (Not every connected algebraic zip datum has that property, for instance, if $L$ and $M$ have root system $A_1$ associated to long and short roots, respectively, and the square of the ratio of the root lengths is different from the characteristic of~$k$.)

Fix a connected reductive algebraic group $G$ over~$k$ and an isogeny $\phi\colon G\to G$. Choose a Borel subgroup $B\subset G$ and a maximal torus $T\subset B$, and let $W$ be the corresponding Weyl group of $G$ and $S$ its set of simple reflections. Choose an element $\gamma\in G$ such that $\phi(\leftexp{\gamma}{B})=B$ and $\phi(\leftexp{\gamma}{T})=T$. Then $\phi\circ\inn(\gamma)\colon \Norm_G(T)\to\Norm_G(T)$ induces an isomorphism of Coxeter systems
\[
\bar\phi\colon (W,S) \iso (W,S).
\]
For any subset $I\subset S$ recall from Subsection~\ref{Weyl} that $P_I$ denotes the standard parabolic of type~$I$. 
Thus the choices imply that $\phi(\leftexp{\gamma}{P_I}) = P_{\bar\phi(I)}$. We denote the unipotent radicals of arbitrary parabolics $P$, $Q$, $P'$, $Q'$ by $U$, $V$, $U'$, $V'$, respectively.

Let $\hat G$ be a linear algebraic group over $k$ having identity component~$G$, and let $G^1$ be an arbitrary connected component of~$\hat G$. Choose an element $g_1\in \Norm_{G^1}(B)\cap\Norm_{G^1}(T)$. Then $\inn(g_1)$ induces an automorphism of $G$ that we use to twist~$\phi$. Let $\delta\colon (W,S)\to (W,S)$ be the isomorphism of Coxeter systems induced by $\inn(g_1)$. 
Then for any subset $I\subset S$ we have $\leftexp{g_1}{P_I} = P_{\delta(I)}$.

Fix subsets $I$, $J\subset S$ and an element $x \in \doubleexp{J}{W}{\delta\bar\phi(I)}$ with $J=\leftexp{x}{\delta\bar\phi(I)}$. Set $y\defeq(\delta\bar\phi)^{-1}(x)\in\nobreak W$.

\begin{lemma}   \label{phieq}
  \begin{itemize}
  \item [(a)] $x\Phi_{\delta\bar\phi(I)}=\Phi_J$.
  \item [(b)] $x\Phi_{\delta\bar\phi(I)}^+=\Phi_J^+$.
  \end{itemize}
\end{lemma}

\begin{proof}
Part (a) follows from  $J=\leftexp{x}{\delta\bar\phi(I)}$. By \eqref{eq:wichar1} the fact that $x\in \doubleexp{J}{W}{\delta\bar\phi(I)}\subset W^{\delta\bar\phi(I)}$ implies $x\Phi_{\delta\bar\phi(I)}^+\subset \Phi^+$. Together with (a) this implies (b).
\end{proof}

\begin{construction}\label{IphixZD}
Set $Q\defeq P_J$ and $P\defeq \leftexp{\gamma\dot y}{P_I}$ and let $L$ be the Levi component of $P$ containing $\leftexp{\gamma\dot y}{T}$. Then $\leftexp{g_1}{\phi}(P) = \leftexp{\dot x}{(}\leftexp{g_1}{\phi}(\leftexp{\gamma}{P_I})) = \leftexp{\dot x}{P_{\delta\bar\phi(I)}}$ and $Q=P_J$ have relative position~$x$. Set $M\defeq\leftexp{g_1}{\phi}(L)$; this is a Levi component of $\leftexp{g_1}{\phi(P)}$ containing $\leftexp{g_1}{\phi}(\leftexp{\gamma\dot y}{T}) = \leftexp{g_1}{\phi}(\leftexp{\gamma}{T}) = \leftexp{g_1}{T} = T$. Since the root system of $M$ is $x\Phi_{\delta\bar\phi(I)}$, Lemma \ref{phieq} shows that it is also the Levi component of $Q$ containing $T$. Let $\leftexp{g_1}{\tilde\phi} \colon P/U\to Q/V$ denote the isogeny corresponding to $\inn(g_1)\circ\phi\,|\,L\colon L\to M$. Then we obtain a connected algebraic zip datum $\CZ\defeq(G,P,Q,\leftexp{g_1}{\tilde\phi})$.
\end{construction}

\begin{lemma}\label{allframe}
The triple $(B,T,\gamma\dot y)$ is a frame of~$\CZ$, and the Levi components determined by it are $M\subset Q$ and $L\subset P$.
\end{lemma}

\begin{proof}
The statements about $M$ and $L$ follow from the inclusions $T\subset M$ and $\leftexp{\gamma\dot y}{T}\subset L$.
They also imply that the isogeny $L\to M$ corresponding to $\leftexp{g_1}{\tilde\phi}$ is simply the restriction of $\inn(g_1)\circ\phi$. Conditions (a) and (b) in Definition \ref{FrameDef} assert that $B\subset Q$ and $\leftexp{\gamma\dot y}{B} \subset P$, which hold by the construction of $Q$ and $P$. Condition (d) translates to $\leftexp{g_1}{\phi}\bigl(\leftexp{\gamma\dot y}{T}\bigr) = T$, which was already shown in~\ref{IphixZD}. 

To prove (c) note first that by Lemma \ref{phieq} we have $x\Phi_{\gamma\bar\phi(I)}^+ = \Phi_J^+$ and therefore $\leftexp{\dot x}{B}\cap M = B\cap M$. The definition of $y$ implies that $\leftexp{g_1}{\phi}(\gamma\dot y\gamma^{-1})\in \dot xT$ and hence
$$\leftexp{g_1}{\phi}\bigl(\leftexp{\gamma\dot y}{B}\bigr) = 
\leftexp{\leftexp{g_1}{\phi}(\gamma\dot y\gamma^{-1})\cdot g_1}{\phi}\bigl(\leftexp{\gamma}{B}\bigr) = 
\leftexp{\leftexp{g_1}{\phi}(\gamma\dot y\gamma^{-1})}{B} = 
\leftexp{\dot x}{B}.$$
From this we can deduce that
$$\leftexp{g_1}{\phi}\bigl(\leftexp{\gamma\dot y}{B}\cap L\bigr) 
= \leftexp{g_1}{\phi}\bigl(\leftexp{\gamma\dot y}{B}\bigr) \cap \leftexp{g_1}{\phi}(L) 
= \leftexp{\dot x}{B}\cap M
= B\cap M,$$ 
proving the remaining condition~(c).
\end{proof}

The automorphism $\psi$ defined in~\eqref{DefinePsi} for the algebraic zip datum $\CZ$ is given by
\begin{equation}\label{psitoisogeny}
\psi := \delta \circ \bar\phi \circ \inn(y) = \inn(x) \circ \delta \circ \bar\phi\colon (W_I,I) \iso (W_J,J)
\end{equation}

\begin{definition}\label{Xdef}
Let $X_{I,\varphi,x}$ be the set of all triples $(P',Q',[g'])$ consisting of parabolic subgroups $P'$, $Q'$ of $G$ of type $I$, $J$ and a double coset $[g'] \defeq V'g'\phi(U')\subset G^1$ of an element $g'\in G^1$ such that
$$\relpos(Q',\leftexp{g'}{\!\phi(P')}) = x.$$
\end{definition}

One readily verifies that the condition on the relative position depends only on~$[g']$, and that $\bigl((g,h),(P',Q',[g'])\bigr) \mapsto \bigl(\leftexp{g}{P'},\leftexp{h}{Q'},[hg'\phi(g)^{-1}]\bigr)$ defines a left action of $G\times G$ on $X_{I,\varphi,x}$. We also have a standard base point $(P,Q,[g_1]) \in X_{I,\varphi,x}$.
One can use the definition of $X_{I,\varphi,x}$ to endow it with the structure of an algebraic variety over~$k$, but in the interest of brevity we define that structure using the following isomorphism:

\begin{proposition}\label{Iquot}
There is a natural $G\times G$-equivariant isomorphism 
$$C_\CZ \stackrel\sim\longto X_{I,\varphi,x},\ 
\bigl(gP,hQ,\ell_{h}\circ \leftexp{g_1}{\tilde\phi} \circ\ell_{g}^{-1}\bigr) \mapsto
\bigl(\leftexp{g}{P},\leftexp{h}{Q},[hg_1\phi(g)^{-1}]\bigr).$$
\end{proposition}

\begin{proof} 
In view of Lemma \ref{quot} the assertion is equivalent to saying that the action of $G\times G$ on $X_{I,\varphi,x}$ is transitive and the stabilizer of $(P,Q,[g_1])$ is~$E_\CZ$. The transitivity follows directly from the definition of the action. For the stabilizer note that $\bigl(\leftexp{g}{P},\leftexp{h}{Q},[hg_1\phi(g)^{-1}]\bigr) = (P,Q,[g_1])$ if and only if $g\in P$ and $h\in Q$ and $Vhg_1\phi(g)^{-1}\phi(U) = Vg_1\phi(U)$. Write $g=u\ell$ for $u\in U$, $\ell\in L$ and $h=vm$ for $v\in V$, $m\in M$. Then the last condition is equivalent to $Vmg_1\phi(\ell)^{-1}\phi(U) = Vg_1\phi(U)$, or again to $m\cdot\leftexp{g_1}{\phi}(\ell)^{-1} \in V{\cdot}\,\leftexp{g_1}{\phi}(U) \cap M$. But for any element $v'\cdot\leftexp{g_1}{\phi}(u') = m'\in V{\cdot}\,\leftexp{g_1}{\phi}(U) \cap M$ we have $\leftexp{g_1}{\phi}(u') = v^{\prime-1}m'\in \leftexp{g_1}{\phi}(U) \cap VM$, and since $M$ is also a Levi component of $\leftexp{g_1}{\phi}(P)$, it follows that $\leftexp{g_1}{\phi}(U) \cap VM = \leftexp{g_1}{\phi}(U) \cap V$ and hence $m'=1$. The last condition is therefore equivalent to $m=\leftexp{g_1}{\phi}(\ell)$. Together this shows that the stabilizer is~$E_\CZ$, as desired.
\end{proof}

\begin{lemma}\label{corr}
For any $w\in \leftexp{I}{W}\cup W^J$ the subset $G^w\subset G$ corresponds via Theorem \ref{orbitbijection} (a) and Proposition \ref{Iquot} to the subset
$$X_{I,\varphi,x}^w\ \defeq\ \bigl\{
\bigl(\leftexp{g}{P_I},\leftexp{g\dot w}{P_J},[g\dot w\dot xg_1 b\phi(\gamma g^{-1})]\bigr)
\;\bigm|\; g\in G,\ b\in B \bigr\} \ \subset\ X_{I,\varphi,x},$$
which is a nonsingular variety of dimension $\dim P + \ell(w)$. 
\end{lemma}

\begin{proof}
Since $(B,T,\gamma\dot y)$ is a frame of~$\CZ$ by Lemma~\ref{allframe}, Theorem \ref{Gw4} for $w \in \leftexp{I}{W}$, respectively \eqref{GwDef2} and Lemma \ref{orbBB} for $w\in W^J$, show that $G^{w} = \orb_\CZ(\gamma\dot yB\dot w)$. In other words $G^w$ is the union of the $E_\CZ$-orbits of $\gamma\dot yb\dot w$ for all $b\in B$. But by \eqref{lambda} and \eqref{rho} we have
\begin{eqnarray*}
\lambda\bigl((\gamma\dot yb)^{-1},\dot w\bigr)
&=& \gamma\dot yb\dot w, \qquad \hbox{and} \\
\rho\bigl((\gamma\dot yb)^{-1},\dot w\bigr) 
&=& \!\! \bigl((\gamma\dot yb)^{-1}P,\dot wQ,\ell_{\dot w}\circ\phi\circ\ell_{(\gamma\dot yb)^{-1}}^{-1}\bigr),
\end{eqnarray*}
and so the $E_\CZ$-orbit of the former corresponds to the $\Delta(G)$-orbit of the latter under the correspondence from Theorem~\ref{orbitbijection}. Moreover, under the isomorphism from Proposition \ref{Iquot} the latter corresponds to the triple
$$\bigl(\leftexp{(\gamma\dot yb)^{-1}}{\!P},\leftexp{\dot w}{Q},[\dot wg_1\phi((\gamma\dot yb)^{-1})^{-1}]\bigr).$$
The definitions of $P$ and $Q$ show that $\leftexp{(\gamma\dot yb)^{-1}}{\!P} = \leftexp{b^{-1}}{\!P_I} = P_I$ and $\leftexp{\dot w}{Q} = \leftexp{\dot w}{P_J}$. The definition of $y$ means that $\leftexp{g_1}{\phi}(\gamma\dot y\gamma^{-1}) = \dot xt$ for some $t\in T$; hence
$$\dot wg_1\phi((\gamma\dot yb)^{-1})^{-1}
= \dot wg_1\phi(\gamma\dot yb)
= \dot w\cdot\dot xt\cdot g_1\cdot\phi(\gamma b\gamma^{-1})\cdot\phi(\gamma)
= \dot w\cdot\dot x\cdot g_1\cdot\leftexp{g_1^{-1}}{\!t}\,\phi(\gamma b\gamma^{-1})\cdot\phi(\gamma).$$
Since $\phi(\leftexp{\gamma}{B})=B$, the factor $b'\defeq\leftexp{g_1^{-1}}{\!t}\,\phi(\gamma b\gamma^{-1})$ runs through $B$ while $b$ runs through~$B$. Thus altogether it follows that $G^w$ corresponds to the union of the $\Delta(G)$-orbits of the triples 
$$\bigl(P_I,\leftexp{\dot w}{P_J},[\dot w\dot xg_1 b' \phi(\gamma)]\bigr)$$
for all $b'\in B$. This union is just the set $X_{I,\varphi,x}^w$ in the lemma. The rest follows from Theorems
\ref{GwVarDim}, \ref{Dualzstablepieces1}, and~\ref{orbitbijection}.
\end{proof}

Combining this with Theorems \ref{thm:zstablepieces1} and \ref{thm:closure} and \ref{orbitbijection} we conclude:

\begin{theorem}
\label{XDecompI}
\begin{itemize}
\item[(a)] The $X_{I,\varphi,x}^w$ for all $w\in {}^IW$ form a disjoint decomposition of $X_{I,\varphi,x}$ by nonsingular subvarieties of dimension $\dim P+\ell(w)$.
\item[(b)] For any $w \in {}^IW$ we have
$$\overline{X_{I,\varphi,x}^w} 
= \coprod_{\substack{w^{\prime}\in {}^IW\\ w^{\prime}\preccurlyeq w}}X_{I,\varphi,x}^{w'}.$$
\end{itemize}
\end{theorem}
Analogously, using Theorems \ref{Dualzstablepieces1} and \ref{dualclosure} and \ref{orbitbijection} we obtain:
\begin{theorem}
\label{XDecompJ}
\begin{itemize}
\item[(a)] The $X_{I,\varphi,x}^w$ for all $W^J$ form a disjoint decomposition of $X_{I,\varphi,x}$ by nonsingular subvarieties of dimension $\dim P+\ell(w)$.
\item[(b)] For any $w \in W^J$ we have
$$\overline{X_{I,\varphi,x}^w} 
= \coprod_{\substack{w^{\prime}\in W^J\\ w^{\prime}\preccurlyeq w}}X_{I,\varphi,x}^{w'}.$$
\end{itemize}
\end{theorem}

\subsection{Frobenius}
\label{RelateMW}

Keeping the notations of the preceding subsection, we now assume that $k$ has positive characteristic and that $\varphi\colon G \to G$ is the Frobenius isogeny coming from a model $G_0$ of $G$ over a finite subfield $\BF_q\subset k$ of cardinality~$q$. Then $G_0$ is quasi-split; hence we may, and do, assume that $B$ and $T$ come from subgroups of $G_0$ defined over $\BF_q$ and therefore satisfy $\phi(B)=B$ and $\phi(T)=T$. We can thus take $\gamma\defeq1$.

In this case, our varieties $X_{I,\phi,x}$ coincide with the varieties $Z_I$ used in~\cite{moonwed} to study $F$-zips with additional structures. The isogeny $\leftexp{g_1}{\tilde\phi}$ in the connected algebraic zip datum $\CZ$ then has vanishing differential; hence $\CZ$ is orbitally finite by Proposition~\ref{dphi=0}. Thus by Theorem \ref{thm:representatives} each $G^w$ is a single $E_\CZ$-orbit, and so by Theorem \ref{orbitbijection} and Theorem \ref{XDecompI} we deduce:

\begin{theorem}
\label{MWOrbits}
\begin{itemize}
\item[(a)] If $\phi$ is the Frobenius isogeny associated to a model of $G$ over a finite field, each $X_{I,\varphi,x}^w$ in Theorem \ref{XDecompI} is a single $\Delta(G)$-orbit. In particular the set
\begin{equation*}
  \bigl\{(P_I,\leftexp{\dot w}{P_J}, [\dot w\dot xg_1]) \bigm| w\in {}^IW \bigr\}
\end{equation*}
is a system of representatives for the action of $\Delta(G)$ on $X_{I,\phi,x}$.
\item[(b)] For any $w\in {}^IW$, the closure of the orbit of $(P_I,\leftexp{\dot w}{P_J}, [\dot w\dot xg_1])$ is the union of the orbits of $(P_I,\leftexp{\dot w'}{\!P_J}, [\dot w'\dot xg_1])$ for those $w'\in {}^IW$ satisfying $w'\preccurlyeq w$.
\end{itemize}
\end{theorem}

Theorem \ref{MWOrbits} (a) was proved in \cite{moonwed}, Theorem 3 and (b) answers the question of the closure relation that was left open there.

\subsection{Lusztig's varieties}

Now we apply the results of Subsection \ref{AZDisog} to the special case $\phi = \id$. In this case we can choose $\gamma\defeq1$ and obtain $\bar\phi = \id$. Then our varieties $X_{I,\phi,x}$ coincide with the varieties $Z_{I,x,\delta}$ defined and studied by Lusztig in \cite{Lusztig:parsheaves2}. There he defines a decomposition of $X_{I,\phi,x}$ into a certain family of $\Delta(G)$-invariant subvarieties. In \cite{He:GStablePieces}, He shows how to parametrize this family by the set $W^{\delta(I)}$. We will denote the piece corresponding to $w\in W^{\delta(I)}$ in this parametrization by $\tilde X^w_{I,\phi,x}$. (In \cite{He:GStablePieces}, He denotes $X_{I,\phi,x}$ by $\tilde Z_{I,x,\delta}$ and $\tilde X^w_{I,\phi,x}$ by $\tilde Z^w_{I,x,\delta}$.) We will show that this decomposition is the same as ours from Theorem \ref{XDecompJ} up to a different parametrization.

\begin{lemma} \label{lem:wjbij}
The map $w\mapsto wx$ induces a bijection $W^{J}\isoto W^{\delta(I)}$.
\end{lemma}

\begin{proof}
Take any $w\in W^J$. Using Lemma \ref{phieq} and \eqref{eq:wichar1} we get $wx\Phi_{\delta(I)}^+=w\Phi_J^+\subset \Phi^+$. By \eqref{eq:wichar1} this shows that $wx\in W^{\delta(I)}$. A similar argument shows that $wx^{-1}\in W^J$ for any $w\in W^{\delta(I)}$, which finishes the proof.
\end{proof}

\begin{theorem}\label{Lusztigstablepieces}
For any $w\in W^J$ we have $X^{w}_{I,\phi,x} = \tilde X^{wx}_{I,\phi,x}$.
\end{theorem}

\begin{proof}
The statement makes sense by Lemma \ref{lem:wjbij}. Let $w\in W^J$ and $w'\defeq wx\in W^{\delta(I)}$. In \cite{He:GStablePieces}, Proposition 1.7, He shows that
\begin{equation*}
  \tilde X^{w'}_{I,\phi,x}=\Delta(G)\cdot
  \bigl\{ (P_I,\leftexp{b\dot w'\dot x^{-1}}{P_J},[b\dot w'g_1 b']) \bigm| b,b'\in B \bigr\}.
\end{equation*}
 (In \cite{He:GStablePieces}, it is assumed that $G$ is semi-simple and adjoint. But this assumption is not needed for the proof of Proposition 1.7 in [loc.\ cit.].) By acting on such a point $(P_I,\leftexp{b\dot w'\dot x^{-1}}{P_J},[b\dot w'g_1 b'])$ with $\Delta(b^{-1})$ and using $w=w'x^{-1}$ we get
 \begin{equation*}
   \tilde X^{w'}_{I,\phi,x}=\Delta(G)\cdot 
   \bigl\{( P_I,\leftexp{\dot w}{P_J},[\dot w\dot x g_1 b']) \bigm| b'\in B \bigr\}.
 \end{equation*}
Since $\gamma=1$, comparison with Lemma \ref{corr} proves the claim.
\end{proof}
From Theorem \ref{XDecompJ} we can now deduce the closure relation between the $\tilde X^w_{I,\phi,x}$:
\begin{theorem}
\label{Lusztigclosure}
  For any $w\in W^{\delta(I)}$ we have
  \begin{equation*}
    \overline{\tilde X^w_{I,\phi,x}}= \coprod_{\substack{w^{\prime}\in W^{\delta(I)}\\ w^{\prime}x^{-1}\preccurlyeq wx^{-1}}}\tilde X^{w'}_{I,\phi,x}.
  \end{equation*}
\end{theorem}
In the special case $x=1$ this result is due to He (see \cite{He:GStablePieces}, Proposition 4.6).

\bibliography{references}
\bibliographystyle{alphanum}

\end{document}